\documentstyle[12pt,amstex,amssymb]{article}

\def\bc{\begin{center}}
\def\ec{\end{center}}
\def\be{\begin{equation}}
\def\ee{\end{equation}}
\def\ben{\begin{enumerate}}
\def\een{\end{enumerate}}
\def\bfg{\begin{figure}}
\def\efg{\end{figure}}
\def\bq{\begin{quote}}
\def\eq{\end{quote}}
\def\bd{\begin{description}}
\def\ed{\end{description}}
\def\this{i.\ e.\ } 

\def\h{\hbar}
\def\p{\partial}
\def\w{\wedge}

\def\dim{\operatorname{dim}}

\def\det{\operatorname{det}}

\def\deg{\operatorname{deg}}
\def\ev{\operatorname{ev}}
\def\ft{\operatorname{ft}}
\def\ct{\operatorname{ct}}
\def\Res{\operatorname{Res}}

\newcommand{\CC}{{\Bbb C}}

\newcommand{\ZZ}{{\Bbb Z}}
\newcommand{\QQ}{{\Bbb Q}}

\newcommand{\lan}{\langle}
\newcommand{\ran}{\rangle}

\renewcommand{\a}{\alpha}
\renewcommand{\b}{\beta}
\renewcommand{\c}{\gamma}
\renewcommand{\d}{\delta}
\newcommand{\e}{\varepsilon}
\newcommand{\f}{\varphi}
\newcommand{\n}{\eta}

\renewcommand{\l}{\lambda}
\newcommand{\m}{\mu}
\renewcommand{\o}{\omega}

\renewcommand{\t}{\tau}

\renewcommand{\v}{\nu}
\newcommand{\x}{\chi}
\newcommand{\xsi}{\xi}

\newcommand{\C}{\Gamma}
\newcommand{\D}{\Delta}
\renewcommand{\L}{\Lambda}

\renewcommand{\S}{\Sigma}

\newcommand{\calo}{{\cal O}}

\newcommand{\calt}{{\cal T}}

\newcommand{\M}{\overline{\cal M}}

\title{Elliptic Gromov-Witten invariants \\ and 
the generalized mirror conjecture.}

\author{Alexander Givental \thanks{ Research supported by NSF grants
DMS-93-21915 and DMS-97-04774} \\ UC Berkeley}

\date{March 12, 1998 
\abstract{\footnotesize 
A conjecture expressing genus $1$ Gromov-Witten invariants in 
mirror-theoretic terms of semi-simple Frobenius structures and
complex oscillating integrals is formulated. The proof of the 
conjecture is given for torus-equivariant Gromov - Witten invariants
of compact K\"ahler manifolds with isolated fixed points and for
concave bundle spaces over such manifolds. Several results on
genus 0 Gromov - Witten theory include: a non-linear Serre duality
theorem, its application to the genus 0 mirror conjecture, 
a mirror theorem for concave bundle spaces over toric manifolds
generalizing a recent result of B. Lian, K. Liu and S.-T. Yau.
We also establish a correspondence (see the extensive
footnote in section 4) between their new proof of the genus 0 
mirror conjecture for quintic 3-folds and our proof of
the same conjecture given two years ago.}
}

\begin{document}

\maketitle

\newpage

\section*{Introduction}

\vspace{0.5in}

Gromov-Witten invariants of a compact symplectic manifold $X$ are defined by
means of enumeration of compact pseudo-holomorphic curves in $X$. 
For any cycle $M$ in the moduli space of genus $g$ Riemann
surfaces with $n$ marked points and for any $n$ cycles in $X$ one can define
a GW-invariant counting those genus $g$ marked pseudo-holomorphic curves in 
$X$ which pass by the marked points through the given cycles in $X$ and whose
holomorphic type belongs to $M$. The handful of GW-invariants thus introduced 
obeys various universal identities which originate from topology of moduli
spaces of Riemann surfaces and constitute a remarkable and fairly sophisticated
algebraic structure. In this paper, we study the structure formed by rational
and elliptic GW-invariants.

The structure of rational GW-invariants alone is well-understood and has been
formalized by B. Dubrovin \cite{Db} in the concept of Frobenius manifolds.
The genus $0$ GW-invariants define on the total cohomology
space $H:=H^*(X)$ a Frobenius manifold structure; roughly speaking, it
consists of the associative commutative
\footnote{ Commutativity and symmetricity should be understood in the sense
of super-algebra since the cohomology space is $\ZZ_2$-graded.} 
{\em quantum cup-product} 
on the tangent spaces $T_tH$ which
is a deformation of the ordinary cup-product, is symmetric with respect 
to the Poincare intersection form $\lan \cdot ,\cdot \ran $ and depends on 
the application point $t\in H$ in such a way that certain integrability 
conditions are satisfied. 

The following observation is a foundation for the so called {\em mirror
conjecture} for Calabi-Yau manifolds and its generalization to arbitrary
symplectic manifolds suggested in \cite{Gi2}: Frobenius manifolds occur
in the fields of mathematics quite remote from symplectic topology or
enumerative geometry, and in particular --- in singularity theory of 
isolated critical point of holomorphic  functions. We outline below the
{\em singularity theory -- symplectic topology} dictionary.

\twocolumn

\begin{enumerate}
\item
A germ $f:(\CC ^m,0) \to (\CC ,0)$ of holomorphic function at an isolated  
critical point of multiplicity $\m$. 

\item
Local algebra $\CC [[z]]/(\p f/\p z)$.

\item
Residue pairing $\lan \phi , \psi \ran :=$
\newline
 $ \int_{ |\p f/\p z_i| =\e_i } 
\frac{\phi (z)\psi (z) dz_1\w ... \w dz_m}
{(\p f/\p z_1) ...(\p f/\p z_m)} $

\item
Parameter space $\L $ of a miniversal deformation 
$f_{\l}(z)$ of the critical point; can be taken in the form
$f_{\l}(z):=$ \newline $f(z)+\l_1\phi_1+...+\l_{\m}\phi_{\m}$
where $\phi_{\a}$ are to represent a basis in the local algebra.

\item
Lagrangian submanifold \newline
 $L\subset T^*\L$ generated by $f_{\l}$, $L:=$ \newline
$ \{ (p,\l)|\exists z: d_zf_{\l}=0, p=\frac{\p f_{\l}}{\p \l} \}$

\item
Critical values $u_1,...,u_{\m}$ of Morse functions $f_{\l}$
at the critical points.

\item
Residue metric $\lan \p_{\l_{\a}}, \p_{\l_{\b}} \ran _{\l} :=$
\newline
$ \sum_{z\in crit (f_{\l})} \frac{(\p_{\l_{\a}} f_{\l})\ 
(\p_{\l_{\b}} f_{\l})}{\det (\p^2 f_{\l}/\p z_i \p z_j)}|_z $ \newline
diagonalized in the basis of non - degenerate critical points 
of a Morse function $f_{\l}$.

\item
Hessians $\D :=\det (\frac{\p^2 f_{\l}}{\p z_i \p z_j})$ 
at the critical points.

\end{enumerate}

\pagebreak

\begin{enumerate}

\item
A compact symplectic manifold $X$. 

\item  Cohomology algebra $H^*(X)$ 

\item  Poincare pairing  on $H^*(X)$, \newline
$\lan \phi , \psi \ran :=\int_X \phi \w \psi $.

\item
The space $H:=H^*(X)$  \newline 
considered as a manifold. 

\item
Spectral variety $L\subset T^*H$ 
of the quantum cup-product $\circ_t$, \newline 
$L\cap T^*_tH := Specm (T_tH, \circ_t )$. \newline
 It is Lagrangian due to the integrability condition mentioned above.

\item
 Function $u:L\to \CC $ such that $du = pdt | L$ considered as a 
multiple - valued function on $H$.

\item
Poincare metric \newline
$\lan \phi ,\psi \ran_t =
\sum_{p\in L\cap T^*_tH} \frac{\phi (p) \psi (p)}{\D (p)}$ \newline
 diagonalized in the basis of \newline
 idempotents of the quantum \newline 
cup-product $\circ_t$  at semi-simple points $t$.  

\item
The function $\D: L\to \CC $ \newline
{\em (quantum Euler class)} representing on $L \subset L\times_H L $
 the cohomology class Poincare-dual  to the diagonal in $X\times X$.  
\end{enumerate}

\onecolumn

The objective of the present paper consists in extending 
the dictionary to include elliptic GW-invariants.

The $3$-valent tensor $\lan a \circ_t b, c \ran $ of structural constants of 
the quantum cup-product on $T_tH$ is actually defined by the formal series:
\[ \lan a \circ_t b, c \ran := \sum_{n=0}^{\infty} (a,b,c,t,...,t)/n! \]
where the GW-invariant $(a,b,c,t,...,t)$ counts the number of rational
curves in $X$ with $n+3$ marked points situated on generic cycles
whose homology classes are Poincare-dual respectively to $a,b,c,t,...,t$.

The genus $1$ GW-invariants in question can be similarly organized into
a uni-valent tensor --- an exact differential $1$-form $dG$ on $H$. The
value of this $1$-form on a tangent vector $a\in T_tH=H^*(X)$ is defined by 
the formal series
\[ i_a dG := \sum _{n=0}^{\infty} [ a, t,...,t ]/n!  \]
where the GW-invariant $[a,t,...,t]$ counts the number of elliptic curves
in $X$ with $n+1$ marked points situated on generic cycles representing
respectively $a,t,...,t$. 

We propose the following construction for the singularity theory counterpart
of the differential $1$-form $dG$ in Gromov-Witten theory.
The critical values $u_1,...,u_{\m}$
of the functions $f_{\l}$ can be taken on the role of local coordinates
on the parameter space $\L $ of the miniversal deformation in the complement
to the {\em caustic} --- the critical value locus of the lagrangian map 
$L\subset T^*\L \to \L$. Consider the {\em complex oscillating integral}
\[ I:=\int_{\C} e^{f_{\l}(z)/\h } v (z,\l) dz_1 \w ... \w dz_m .\]
The partial derivative $\h \p I/\p u_{\a}$ of the complex
oscillating integral can be expanded into the stationary phase asymptotical 
series
\[ \h^{m/2} e^{u_{\a}/\h}\frac{(\p f_{\l}/\p u_{\a}) |_{z_{crit}}}
{\sqrt{\D_{\a}}}\ (1+ \h R_{\a} + o(\h) ) \]
near the non-degenerate critical point $z_{crit}$ corresponding to 
{\em the same} critical value $u_{\a}$. The {\em asymptotical coefficients} 
$R_{\a}$ actually
depend only on the $4$-jet of $f_{\l}$ and the $2$-jet of $v $ at 
$z_{crit}$. In terms of the critical values 
$u_{\a}$, the Hessians $\D_{\a}$ and the asymptotical coefficients 
$R_{\a}$ the differential
$1$-form $dG$ is described by the formula:
\[  (*) \ \ \ dG = \sum _{\a =1}^{\m} 
(\frac{1}{48} d\log \D_{\a} + \frac{1}{2}R_{\a}du_{\a}) .\] 
 
Our proposal has implications in both singularity and Gromov-Witten theory.

In singularity theory, the residue metric on $\L$ 
(which is the counterpart of the flat Poincare metric on $H$) has no reason
to be flat. However, according to K. Saito  theory of primitive 
forms \cite{Sa} one can choose the holomorphic volume form 
$v(z,\l) dz_1\w ... \w dz_m$ (called {\em primitive})
in such a way that the corresponding 
residue metric is flat. Moreover, Saito's theory can be reformulated as the
theorem that the above dictionary introduces a Frobenius structure on $\L $
{\em provided that} $(z_1,...,z_m)$ everywhere in the dictionary means a 
unimodular coordinate system with 
respect to the primitive volume form.  The same primitive form should be used
in the definition of complex
oscillating integrals involved into our construction of the $1$-form $(*)$.
With this hypothesis in force, we arrive to the following

\medskip

{\bf Conjecture 0.1.} {\em The $1$-form (*) satisfies all axioms for the genus
$1$ GW-invariant $dG$. In particular, E. Getzler's relation \cite{Ge} holods
true for $(*)$.}

\medskip

{\em Remark on examples.} 
We will see in Section $1$ from the theory of Frobenius 
structures that differentials of the asymptotical coefficients $R_{\a}$ are 
expressible via the Hessians $\D_{\b}$ by
\[ dR_{\a} =\frac{1}{4}\sum_{\b} (\p_{\a} \log \D_{\b})(\p_{\b} \log \D_{\a})
 (du_{\b}-du_{\a}) ,\]
where $\p_{\c}$ means partial derivatives in the coordinate system 
$(u_1,...,u_{\mu})$. This allows to compute $R_{\a}$ if the Hessians are known
as functions of all $u_{\c}$. However in applications to Gromov -- Witten 
theory $\D_{\b}$ are usually known only along some subspace in $H$, 
and the asymptotical coefficients are to be computed independently.  
   All examples considered in this paper are elementary and have $\mu =2$. In
such a case there are two coordinates $u_{\c}$ which we usually denote 
$u_{\pm}$. The functions $R_{\pm}$ and $\D_{\pm}$ depend only on the
difference $u=u_{+}-u_{-}$. In the most examples we will have 
$R_{\pm}=\pm R (u)$ and $\D_{\pm}=\pm D (u)$. Then the formula $(*)$ reduces to
\[ dG = \frac{1}{24} d\log \D (u) + \frac{1}{2} R(u)du, \ \text{where} \ 
\frac{dR}{du}=\frac{1}{4} (\frac{d \log \D (u)}{d u})^2 .\]
We use this method in the following example, but in some other examples 
we will present alternative techniques for computing asymptotical coefficients 
in order to illustrate computational tools available in applications to 
sympectic topology.

\medskip

{\em Example.} The critical point $f=x^3$ of type $A_2$ has the miniversal
deformation $x^3-t_1x+t_0$. At the critical point $x_{\pm}=\pm (t_1/3)^{1/2}$
we have the critical value $u_{\pm}=t_0-2x_{\pm}^3$ and the Hessian 
$\D_{\pm}=6x_{\pm}$. Thus $\D (u)= (-u/4)^{1/3}$, $d (\log \D )/du=1/3u$,
$R(u)=\int du /(36 u^2) =-1/(36 u) + const $ where $const =0$ by 
quasi-homogeneity. Therefore
\[ dG = \frac{1}{24}\frac{du}{3u} - \frac{1}{2}\frac{du}{36u} =0 .\]
This result implies (by Hartogs principle) that for any isolated critical 
point $dG$ defined by $(*)$ outside the caustic in the base of miniversal
deformation extends holomorphicly to the whole base, and that $dG=0$ for 
all simple singularities $A_{\mu}, D_{\mu}, E_{\mu}$ (since $dG$ has zero 
quasi-homogeneity degree). The last conclusion agrees with E. Getzler's
relation.
\footnote{I am thankful to E. Getzler for correcting a mistake in my original
computation.}    

\medskip

In Gromov-Witten theory, the counterpart of complex oscillating integrals of 
singularity theory can be defined, as we shell see in the next section, 
entirely in terms of genus $0$ GW-invariants. In particular, the coefficients
$R_{\a}$ and the $1$-form $(*)$ can be defined in intrincic terms of the 
Frobenius structure on $H$ {\em provided} that the algebras $(T_tH, \circ_t)$
are semi-simple for generic $t\in H$.

\medskip

{\bf Conjecture 0.2.} {\em The elliptic GW-invariant $dG$ of a compact 
symplectic manifold $X$ with genericly semi-simple quantum cup-product
is expressed by the formula $(*)$ in terms of rational GW-invariants.} 

\medskip

In the rest of the paper we present our evidence in favor of the conjectures.

In Section $1$ we review some definitions and results of genus $0$
GW-theory, give a more precise formulation of Conjecture $0.2$ and verify 
it directly in the example $X=\CC P^1$.

In Section $2$ we generalize the conjecture to the case of equivariant
GW-theory on a K\"ahler manifold $X$ provided with a Killing torus action with
isolated fixed points only (toric and flag manifolds are main examples).
The corresponding Frobenius structure is genericly semi-simple.

In Section $3$ we prove the equivariant version of the conjecture.

In Section $4$ we introduce and study GW-invariants of the non-compact 
manifolds which are total spaces of sums of negative line bundles over
toric manifolds. Results of Sections $2$ and $3$ extend easily to such spaces.
The key new point in this Section is the {\em mirror theorem} saying that
those rational GW-invariants of such bundles which play the role of
oscillating integrals in the Frobenius structure {\em are equal} 
to certain oscillating integrals defined in the spirit of toric
hyper-geometric functions. In particular, the version of Conjecture
$0.1$ for such oscillating integrals holds true.

Results of Section $4$ provide a new illustration to the so called
mirror phenomenon: not only the GW-invariants of a manifold $X$ form
a structure {\em analogous} to the one observed in singularity theory, but 
the GW-invariants are explicitly {\em expressed} in terms of a specific 
datum of singularity theory type called the {\em mirror partner} of $X$.
      
In Section $5$ we deal with toric {\em super}-manifolds which are objects
dual to the bundles of Section $4$ and whose GW-invariants are to 
coincide with GW-invariants of toric complete intersections. 
We invoke the nonlinear Serre duality theorem \cite{Gi1} 
(which relates genus $0$ GW-invariants of a super-manifold and of the dual 
bundle space) in order to give a new proof of the mirror theorem 
for toric complete intersections \cite{Gi3}. 
The role of hyper-geometric functions is
more transparent in this version of the proof. We believe that elliptic
GW-invariants of toric complete intersections are expressible in terms
of the rational GW-invariants of the corresponding equivariant super-manifolds.
However we are not ready to report on such applications because of the 
difficulties we explain in the end of Section $5$. 

\medskip

I am thankful to organizers and participants of the
summer-97 Taniguchi Symposium where the main results of this paper were first
announced.

\section{Gromov-Witten invariants and semi-simple Frobenius structures}

We review here some basic properties of Gromov-Witten invariants of
compact symplectic manifolds \cite{KM, BM, BF, FO, LT, R1, Db}. 

\medskip

{\bf Stable maps.} 
Let $(\S ,\e)$ denote a {\em prestable marked curve}, that is
a compact connected complex curve $\S $ with at most
double singular points and an ordered $n$-tuple $(\e_1,...,\e_n)$ of
distinct non-singular marked points. The {\em genus} of $(\S,\e)$ is defined
as $g=\dim H^1(\S,\calo_{\S})$. The {\em degree} of a holomorphic map
$f: (\S,\e)\to X$ to a compact (almost) K\"ahler manifold $X$ is defined as
the total homology class $d\in H_2(X,\ZZ )$ the map $f$ represents. 
Two maps $f: (\S,\e)\to X$ and $f': (\S',\e')\to X$ are called
{\em equivalent} if there exist an isomorphism $\f: (\S ,\e)\to (\S',\e')$
such that $f=f'\circ \f$.
A holomorphic map $f: (\S ,\e)\to X$ is called {\em stable} if it has no
non-trivial infinitesimal automorphisms. The set of equivalence classes of
degree $d$ stable holomorphic maps to $X$ of genus $g$ curves with $n$ marked 
points is denoted $X_{g,n,d}$ ( and called {\em 
moduli space of stable maps} ). According to \cite{BF, FO, LT, R1} the moduli
spaces have a natural structure of compact orbi-spaces, complex-analytic if
$X$ is K\"ahler. If $X=pt$, the spaces $X_{g,n,0}$ coincide with the 
Deligne-Mumford compactifications $\M_{g,n}$ of moduli spaces of marked 
Riemann surfaces and are orbifolds of dimension $3g-3+n$ (unless empty, which
happens for $g=0, n<3$ and $g=1, n=0$). For any $X$ degree $0$ stable maps 
form the moduli spaces $X_{g,n,0}=X\times \M_{g,n}$.  
  
One introduces the following tautological maps:

- {\em evaluation maps} $\ev=(\ev_1,...,\ev_n): X_{g,n,d}\to X^n$
defined by evaluating stable maps at the marked points;

- {\em forgetting maps} $\ft_i: X_{g,n+1,d}\to X_{g,n,d},\ i=1,...,n$,
well-defined (unless $d=0$ and $\M_{g,n}$ is empty) by
forgetting the marked point $\e_i$ followed by contracting those irreducible
components of $\S$ which have become unstable;

- {\em contraction maps} $\ct : X_{g,n,d}\to \M_{g,n}$ defined by forgetting
the map $f: (\S,\e)\to X$ followed by 
contracting unstable irreducible components of the marked curve $(\S,\e)$.

The diagram formed by the forgetting map $\ft_{n+1}: X_{g,n+1,d}\to X_{g,n,d}$
and by the evaluation map $\ev_{n+1}: X_{g,n+1,d}\to X$ is called the
{\em universal stable map}: the fibre of $\ft_{n+1}$ over the point represented
by a stable map $f: (\S,\e)\to X$ is canonically identified with (the quotient
of) the curve $\S $ (by the discrete group $Aut (f)$ of automorphisms of the 
map $f$ if this group is non-trivial), and the restriction of $\ev_{n+1}$ to
the fibre (lifted to $(\S,\e)$) is equivalent to $f$. In particular, the
sections $\e_1,...,\e_n: X_{g,n,d}\to X_{g,n+1,d}$ defined by the marked
points play the role of {\em universal} marked points on the universal stable
map. 

One introduces the {\em universal cotangent line} $l_i$ which is a 
line (orbi-)bundle over $X_{g,n,d}$ with the fibre $T^*_{\e_i}\S $ at the
point $[f]$ and defined as the conormal bundle to the universal marked point 
$\e_i$. The $1$-st Chern classes $c^{(1)},...,c^{(n)}$ of the orbi-bundles
$l_1,...,l_n$ are well defined over $\QQ $.

\medskip

{\bf Gromov -- Witten invariants.}
Let $T(c)=t^{(0)}+t^{(1)}c+t^{(2)}c^2+...$ denote a formal power series
with coefficients $t^{(i)}$ in the cohomology algebra 
\footnote{We always assume rational coefficients unless otherwise specified
explicitly.} $H^*(X)$. Given $n$ such series $T_1,...,T_n$, one introduces 
the genus $0$ Gromov-Witten invariant of $X$ by
\[ (T_1,...,T_n)_d := \int_{[X_{0,n,d}]}
(\ev_1^* T_1)(c^{(1)})\w ... \w (\ev_n^* T_n)(c^{(n)}) .\]
Here integration means evaluation of a cohomology class on the {\em virtual}
fundamental class of the moduli space. If $X$ is a {\em convex} K\"hler 
manifold, \this if the tangent bundle $\calt_X$ is spanned by global vector 
fields on $X$, then the genus $0$ moduli spaces $X_{0,n,d}$ are known to be
compact complex orbifolds of complex dimension 
$\lan c_1(\calt_X), d \ran +\dim_{\CC} X +n-3$ 
(see \cite{BM}), and $[X_{0,n,d}]$ is the
findamental class of the orbifold which is well-defined over $\QQ $. In general
the moduli spaces can have many irreducible components of different dimensions
with nasty singularities. Newertheless one can endow them with rational 
virtual fundamental classes of {\em Riemann-Roch} dimension 
\[ \dim_{\CC} [X_{g,n,d}] = \lan c_1(\calt_X), d\ran +(1-g)(\dim_{\CC} X-3) 
+ n \]
in such a way that the axioms \cite{KM} of Gromov-Witten theory are satisfied.
We refer the reader to \cite{BF,FO,LT,R1} for several constructions of the
virtual fundamental classes an for their properties. Using the classes
$[X_{g,n,d}]$ one can introduce higher genus GW-invariants.
\footnote{More general GW-invariants (like $A (T_1,...,T_n)_d$) corresponding
to a choice of a cohomology class $A\in H^*(\M_{g,n})$ are defined by adding
the factor $\ct^*A $ to the integrand.} 
In this paper we will use the notation
\footnote{We will add the super-script indicating the number of arguments 
as in $(T,...,T)_d^n$ or $[T,...,T]_d^n$ when the number would otherwise be
ambiguous.}  
$[T_1,...,T_n]_d$ for the genus $1$ GW-invariants of
$X$. In the case when the series $T_i=t_i\in H^*(X)$ do not depend on $c$
the GW-invariants $(t_1,...,t_n)_d$ (resp. $[t_1,...,t_n]_d$) have the 
enumerative meaning of the number of degree $d$ rational (resp. elliptic) 
curves in $X$ passing through $n$ generic cycles Poincare-dual to 
$t_1,...,t_n$.

\medskip

We are going to use several universal identities between GW-invariants.
\footnote{While the identities are frequently used and their origin is 
well-known and explained for instance in \cite{Gi1}, the actual proofs
depend on details of the definition of the virtual fundamental cycles.
A definition sufficient for our purposes is contained in \cite{LT}.
It is based on the observation that the standard in algebraic geometry
construction of the normal cone to
the zero locus $Z$ of an algebraic section of a vector bundle is intrincic 
with respect to the following data: (the variety $Z$, the complex $E\to F$
of vector bundles over $Z$ defined by the linearization of the section).
The kernel $T=ker(E\to F)$ is the algebraic tangent space to $Z$, and the 
cokernel $N$ is called the obstruction space. The construction is adjusted to
the orbi-bundle setting by applying it equivariantly on the total space of 
a suitable principal orbi-bundle. In the case when $Z$ is a moduli space
of stable maps the tangent and obstruction spaces are already defined
(roughly, as the kernel and cokernel of the Cauchy - Riemann operator).
One can explicitly point out a global resolution $T\to E\to F\to N$
of $T$ and $N$ by a suitable complex $E\to F$ of orbi-bundles. This defines
the intrincic {\em normal subcone} in $F$, and the virtual fundamental cycle is
defined as the intersection of this subcone with generic sections of $F$.
As it is stated in \cite{LT}, with this definition the standard arguments 
justifying the axioms \cite{KM} of GW-theory (they include the string,
divisor, and WDVV-equations) go through. We do not know however a convenient
reference where the details are written down.}

The {\em WDVV-equation} says that the following sum is totally symmetric
in $A, B, C, D$:
\[  \sum_{n'+n''=n} \frac{1}{n'! n''!} \sum_{d'+d''=d} \sum_{\v\v'}
(A,B,T,...,T,\phi_{\a})_{d'}^{n'+3} \n^{\a\b} 
(\phi_{\b}, T,...,T,C,D)_{d''}^{n''+3}    .\]
Here $\{ \phi_{\a}\} $ is a basis in $H^*(X)$, and 
$\sum_{\a\b}\n^{\a\b} \phi_{\a} \otimes \phi_{\b} $ represents the class
in $H^*(X\times X)$ Poincare-dual to the diagonal. In particular the matrix
$(\n^{\a\b})$ is inverse to the intersection matrix
\[ \n_{\a\b}:=\lan \phi_{\a},\phi_{\b} \ran := \int_X \phi_{\a}\w \phi_{\b} .\]

The {\em string} and {\em divisor} equations read respectively:
\[ (1,T_1,...,T_n)_d = \sum _{i=1}^n (T_1,...,DT_i,...,T_n)_d ,\]
\[ (p,T_1,...,T_n)_d = \sum_{i=1}^n (T_1,...,p DT_i,...,T_n)_d+
\lan p , d\ran (T_1,...,T_n)_d \]
where $n$ should be at least $3$ if $d=0$, $DT$ denotes the series 
$(T(c)-T(0))/c $, and $p\in H^2(X)$. The string and divisor
equations hold true for genus $1$ GW-invariants $[ ... ]_d$ (with $n\geq 1$ if
$d=0$) and for GW-invariants of higher genus as well.

\medskip

{\bf Gromov -- Witten potentials.} 
The WDVV-, string and divisor equations have several important interpretations
in terms of the following generating functions for genus $0$ GW-invariants.
The genus $0$ {\em GW-potential} is defined as the formal function of 
$t\in H^*(X)$:   
\[ F(t,q)=\sum_{n=0}^{\infty} \frac{1}{n!} \sum_{d\in \L} q^d (t,...,t)_d^n .\]
Here $\L$ denotes the {\em Mori cone} of $X$, the semi-group in the lattice
$H_2(X,\ZZ )$ generated by those degrees of holomorphic curves in $X$ for
which the virtual fundamental classes $[X_{g,n,d}]$ are non-zero. All $d\in \L$
have non-negative coordinates $(d_1,...,d_r)$ with respect to a suitable basis
\footnote{We mod out the torsion in $H_2(X,\ZZ)$ and thus treat it as a free
abelian subgroup in $H_2(X,\QQ)$ of rank $r$.}
$p_1,...,p_r\in H^2(X,\ZZ )$. Thus the formal completion $\QQ [[\L ]]$ of the 
semigroup algebra of $\L $ is identified with a subalgebra in 
$\QQ [[q_1,...,q_r]]$. The symbol $q^d=q_1^{d_1}...q_r^{d_r}$ stands therefore
for the element $d\in \L$ in the semi-group algebra.

Denote $t_0$ the coordinate on $H^0(X)$, $(t_1,...,t_r)$ --- the coordinates
on $H^2(X)$ with respect to the basis $\{ p_i \}$, so that 
$t:=\sum t_{\a} \phi_{\a} = t_0+t_1p_1+...+t_rp_r+...$, 
and assume that the basis $\{ \phi_{\a} \} = (1, p_1, ..., p_r, ...) $ 
is graded. 

The GW-potential $F$ has the following obvious properties:

- $F$ is homogeneous of degree $3-\dim _{\CC} X$ with respect to the 
grading
\[ \deg t_{\a} = 1-\deg \phi_{\a} /2, \ \deg q^d = \lan c_1(\calt_X), d\ran ,\]

- $F(t,0)=\int _X t\w t\w t \ /6$,

- $\tilde{F}(t,q):=F(t,q)-F(t,0)$ does not depend on $t_0$ and satisfies
$q_i\p \tilde{F} /\p q_i = \p \tilde{F} /\p t_i$ (string and divisor equations
for $(t,...,t)_d$) and thus  
\[ \tilde{F} (t+\sum \t_i p_i, q) = \tilde{F} (t, q\exp \t ) .\] 

We will also make use of the following generating functions:
\[ S_{\a\b}(t,q,\h):=\sum _{n=0}^{\infty} \frac{1}{n!} \sum_{d\in\L} q^d
(\phi_{\a}, t,...,t,\frac{\phi_{\b}}{\h-c})_d^{n+2} ,\]
\[ V_{\a\b}(t,q,x,y):=\sum_{n=0}^{\infty} \frac{1}{n!}\sum_{d\in\L} q^d
(\frac{\phi_{\a}}{x-c},t,...,t,\frac{\phi_{\b}}{y-c})_d^{n+2} .\]
The ill-defined terms in these series are to be replaced as follows:
\[ (\phi_{\a},\frac{\phi_{\b}}{\h-c})_0:= \n_{\a\b},\ 
   (\frac{\phi_{\a}}{x-c},\frac{\phi_{\b}}{y-c})_0 := \frac{\n_{\a\b}}{x+y} .\]
The tensor fields
\[ \sum_{\e\e'}S_{\e\e'} dt_{\e}dt_{\e'} \ \text{and}\ 
\sum_{\e\e'} V_{\e\e'} dt_{\e}dt_{\e'} \] 
have degrees respectively $2-\dim_{\CC} X$ and $1-\dim_{\CC} X$ 
with respect to the above grading and $\deg \h =\deg x=\deg y=1$.

\medskip

In the following description of some identities between the GW-potentials 
$F, S, V$ we will denote $\p _{\a}$ the partial derivatives 
$\p/\p t_{\a}$  with respect to a basis $\{\phi_{\a} \}$ in $H^*(X)$. 
In the formulas below we will ignore
the signs which may occur due to $\ZZ_2$-grading in cohomology and therefore
assume that $H^*(X)$ has no odd part. 

\medskip

(1) Put $F_{\a\b\c}:=\p_{\a}\p_{\b}\p_{\c} F_0 (t,q)$. {\em The} 
WDVV-equation for the GW-potential $F$ reads:
\[ \sum_{\e\e'} F_{\a\b\e}\n^{\e\e'}F_{\e'\c\d} \ \text{is symmetric
in} \ \a,\b,\c,\d \ .\]
This identity is interpreted as associativity of the {\em quantum cup-product}
$\circ : H^*(X)\otimes H^*(X) \to H^*(X)$ defined by the structural constants
\[ \lan \phi_{\a}\circ \phi_{\b}, \phi_{\c} \ran := F_{\a\b\c} \]
(depending on parameters $t=\sum t_{\a}\phi_{\a}$ and $q=(q_1,...,q_r)$).
The quantum cup-product is commutative, symmetric relatively the 
intersection form,
\[ \lan \phi_{\a}\circ \phi_{\b},\phi_{\c}\ran =\lan \phi_{\a},\phi_{\b}\circ
\phi_{\c}\ran \ ,\]
and the unity $1\in H^*(X)$ remains the unity for the quantum cup-product.

(2) The WDVV-equation for $F$ is 
also interpreted as integrability of the following system of linear PDE for a 
vector-function of $t$ (depending also on the parameters $q$
and $\h$) with values in the cohomology space of $X$:
\[ (**) \ \ \h \p_{\a} \vec{s} = \phi_{\a} \circ \vec{s} \ .\]
The formally adjoint system with respect to the intersection form is 
$-\h \p_{\a} \vec{s} =\phi_{\a} \circ \vec{s}$:
\[ \forall \a \ \ \p_{\a} \lan \vec{s} (t,q,\h), \vec{s} (t,q, -\h) \ran =0.\]
 
(3) The generating functions $(S_{\b\c})$ form a fundamental solution
matrix $S$ for the system of PDE:
\[ \h \p_{\a} S_{\b\c} = \sum_{\e\e'} F_{\a\b\e} \n^{\e\e'} S_{\e'\c} \ .\]
Namely,  WDVV-equations imply 
$\phi_{\a}\circ \p_{\b} S= \phi_{\b}\circ \p_{\a} S$
while the string equation implies that $\h \p_0 S = S$.

\medskip

(4) Application of WDVV-equations implies
\[ \p_0 V_{\a\b} (t,q,x,y) = \sum_{\e\e'} \p_0 S_{\e\a}(t,q,x) \n^{\e\e'}
\p_0 S_{\e'\b}(t,q,y) .\]
Together with the string equation this yields the {\em unitarity condition}
\[ \sum_{\e\e'} S_{\e\a}(t,q,\h) \n^{\e\e'} S_{\e'\b}(t,q,-\h)=\n_{\a\b} \]
and the relation
\[ V_{\a\b}(t,q,x,y)=\frac{1}{x+y} \sum_{\e\e'} \n^{\e\e'} S_{\e\a}(t,q,x)
S_{\e'\b}(t,q,y) .\]

(5) The divisor equation with $p=\sum \t_i p_i$ 
applied to $S(t,q,\h)$ shows that
\[ S_{\a\b}(t+p,q,\h)= \sum_{\e} S_{\a\e}(t,q e^{\t},\h) \n^{\e\e'} 
\lan \phi_{\e'} e^{p\t /\h}, \phi_{\b} \ran .\]
This property together with the unitarity condition and the asymptotics
$S_{\a\b}|_{t=0,q=0} = \n_{\a\b}$ uniquely specifies $S$ among fundamental
solutions of the differential system $(**)$.

\medskip

{\bf Frobenius manifolds.} 
A {\em Frobenius algebra} structure on a vector space consists of a 
commutative associative multiplication $\circ $ with unity $1$ and a linear 
function $\a $ such that $\lan u, v\ran := \a (u\circ v)$ is a non-degenerate
bilinear form.

A {\em Frobenius structure} on a manifold $H$ is a field of Frobenius algebra
structures on the tangent spaces $T_tH$ satisfying the following 
integrability conditions:

(a) the metric $\lan \cdot ,\cdot \ran $ is flat: $\nabla ^2=0$,

(b) the unity vector field $1$ is covariantly constant: $\nabla 1 =0$,

(c) the $1$-st order linear PDE system for sections $s$ of $TH$ defined by
$\h \nabla_w s = w\circ s$ is consistent for any $\h \neq 0$.

The Frobenius manifold is said {\em conformal of dimension $D\in \QQ$} 
if it is provided with a vector field $E$ (called {\em Euler}) such 
that the tensor fields  $1, \circ $ and $\lan \cdot , \cdot \ran $ are
eigen-vectors of the Lie derivative operator $L_E$ with the eigen-values
respectively $-1, 1$ and $2-D$. 

In flat coordinates $\{ t_{\a} \}$ of the metric the condition (c) can be
reformulated as flatness for any $\h $ of the connection
\[ (1)\ \ \  \nabla_{\h} := \h d - \sum_{\a} A_{\a} (t) dt_{\a} \w \]
where $A_{\a}$ are the multiplication operators $\p_{\a} \circ_t$.

Using the property of the structural constants $F_{\a\b\c}$ of 
the quantum cup-product on $H^*(X)$ to depend on $q_i, t_i, i=1,...,r$ only
in the combinations $q_i\exp t_i$, we see that the quantum cup-product and 
the Poincare pairing define on $H=H^*(X,\CC )/ 2\pi i H^2(X,\ZZ )$ 
the structure of a (formal) Frobenius manifold of conformal dimension 
$D=\dim_{\CC} X$ with respect to the Euler vector field
\[ E=t_0 \p_0 +\sum_{i=1}^r c_i \p_i +\sum_{\a: \deg t_{\a}<0} 
\deg (t_{\a}) t_{\a} \p_{\a} .\]
Here $\sum c_i p_i$ is the $1$-st Chern class of $\calt_X$.

\medskip

Given a pensil of flat connections $\nabla _{\h}$ one can study asymptotical
behavior of horizontal sections as $\h \to 0$. The
asymptotics is described by the following data.

- The characteristic Lagrangian variety $L\subset T^*H$ defined as 
the spectrum $Spec (Vect(H),\circ )$ of the algebra of vector
fields on $H$ with the multiplication $\circ $.
Flatness of $\nabla_{\h}= \h d - A^1$ is equivalent to $A^1\w A^1=0$
and $dA^1=0$. The first condition means commutativity $[A_{\a},A_{\b}]=0$ while
the second one {\em implies} that $L$ is Lagrangian at generic points 
\cite{GK}.

- The function $u$ on $L$, may be multiple-valued, 
defined as a potential for the action $1$-form 
$\sum p_{\a}dt_{\a} $ on $T^*H$ restricted to $L$. In our case of conformal
frobenius structures $u$ can be chosen as the restriction to $L$ of 
the function $\sum c_i p_i +\sum (\deg t_{\a})p_{\a}t_{\a}$ on $T^*H$ 
defined by the Euler vector field $E$.

- The function $\D$ on $L$ defined by the metric $\lan p, p\ran $ on $T^*H$.
It is the restriction to $L$ of $\sum p_{\a} \n^{\a\b} p_{\b} $.

A point $t\in H$ where the algebra $(T_tH,\circ_t)$ is semi-simple is called
{\em semi-simple}. In a neighborhood of a semi-simple point the characteristic
variety $L$ consists of $N=\dim H$ sections of $T^*H$ which span each $T^*_tH$
so that the corresponding branches $u_{\a}, \a=1,...,N$, of the function $u$ 
form a local coordinate system on $H$ called {\em canonical}. 
The vector fields 
$f_{\a}=\D ^{-1/2} \p /\p u_{\a}$ form an orthonormal basis diagonalising 
$\circ$. Let $\Psi $ denote the transition matrix from the basis 
$\p /\p t_{\a}$ to $f_{\b}$: $\p/\p t_{\a} = \sum \psi_{\a\b} f_{\b}$.

\medskip

{\bf Proposition 1.1.} {\em
In a neighborhood of a semi-simple point there exists a  
fundamental solution of the system  $\nabla_{\h} s=0$ represented 
by the asymptotical series
\[ (2)\ \ \ \Psi (1 + \h R^{(0)} + \h^2 R^{(1)} + ...) \exp (U/\h ) \]
where $U=diag(u_1,...,u_N)$ is the diagonal matrix of canonical coordinates.}

\medskip

{\em Proof.} Indeed, substituting the series into the equation we obtain the 
chain of equations 
\[ A^1 \Psi = \Psi dU ,\ \Psi^{-1} d\Psi = [dU, R^{(0)}] , \] 
\[ D R^{(0)}=[dU, R^{(1)}], \ ..., \ DR^{(k)}=[dU , R^{(k+1)} ], ...\]
where the connection operator $D=d+\Psi^{-1}d\Psi \w $ is flat and anti-
commutes with $dU\w $. 
The first equation means that $\Psi $ diagonalizes $A^1$ to $dU$. Columns of
$\Psi $ form an orthogonal basis since $A^1$ is self-adjoint with respect
to $\lan \cdot ,\cdot \ran$. The next equation requires
the columns to be normalized to constant lengths and expresses off-diagonal
entries of $R^{(0)}$ via $\Psi $. In particular $R^{(0)}$ is symmetric.
The diagonal entries of $R^{(0)}$ can be found by integration from the next
equation: 
\[ dR^{(0)}_{ii}=\sum_l R^{(0)}_{il} (du_l-du_i) R^{(0)}_{li} .\]
Closedness of the RHS is easy to derive directly from the flatness $D^2=0$.
Continuing the inductive procedure, we express the off-diagonal part of 
$R^{(k+1)}$ via $R^{(k)}$ algebraicly from $[dU, R^{(k+1)}]=DR^{(k)}$,
and find the diagonal part of $R^{(k+1)}$ by integration from the next 
equation. Let us check compatibility conditions.

First, $DR^{(k)}$ has the zero diagonal (induction hypothesis) and thus 
is a commutator with $dU$ due to De Rham lemma: the anti-commutator 
\newline $\{ dU, DR^{(k)} \}=D [R^{(k)},dU]=-D^2 R^{(k-1)} =0$. 

It remains to verify exactness of $\sum_l R^{(0)}_{il} (du_l-du_i) 
R^{(k+1)}_{li}$. It can be reformulated as 
$d (R^{(0)} [dU, R^{(k+1)}])_{diag} =0$. We have:
\[ d(R^{(0)} [dU, R^{(k+1)}])=(dR^{(0)})\w DR^{(k)}+ R^{(0)}dDR^{(k)} \]
\[ (dR^{(0)}-[dU,R^{(0)}])\w DR^{(k)}=
(DR^{(0)})^t\w DR^{(k)}=[(R^{(1)})^t, dU]\w [dU, R^{(k+1)}] \]
which has zero diagonal entries. $\square$
  
\medskip 

{\em Remarks.} (1) The asymptotical solution of Proposition $1.1$ is not 
unique. First, the canonical coordinates are defined up to a constant summand.
When the choice has been made, the matrix $\Psi $ of eigen-vectors 
is defined up to the right multiplication by a constant diagonal matrix. 
Such a multiplication conjugates all $R^{(k)}$ by this matrix and thus does 
not change the diagonal entries. Finally, 
another choice of intgration constants
for $R^{(k)}_{ii}$ gives rise to the right multiplication of the whole series
by a diagonal matrix $diag\ (C_1(\h), ...,C_N (\h)) $, 
where $C_i=1+ c_i^{(0)} \h + c_i^{(1)}\h^2+...$ and $c_i^{(j)}$ are constants.

(2) In the canonical coordinate system $(u_1,...,u_N)$ flatness of the
connection $D$ reads: 
\footnote{The last relation means that the vector $\sum \p_l $ represents
the unity in $T_tH$ and holds true for all $R^{(k)}_{ij}$.}
\[ \p_l R^{(0)}_{ij}=R^{(0)}_{il}R^{(0)}_{lj} \ , l\neq i,j,\ 
\sum_l \p_l R^{(0)}_{ij} =0. \  \]
Since $R^{(0)}$ is symmetric, this implies that the $1$-form 
\[ dR:= \sum_{i=1}^N R^{(0)}_{ii} du_i \]
is closed.

(3) The flat metric on $H$ is diagonal in the canonical coordinate system
and takes on $\sum_i \D_i^{-1} (du_i)^2$. The connection $D$ is the Levi-Civita
connection of this metric written in the basis of vector fields 
$\D_i^{1/2} \p /\p u_i$. Thus the coefficients $R_{ij}$ with $i\neq j$ can be
computed in terms of the metric. This observation leads to the formula 
mentioned in the introduction:
\[ dR^{(0)}_{ii}=\frac{1}{4}\sum_j (\p_i \log \D_j) (\p_j \log \D_i) 
(du_j-du_i) .\]   

(4) In the conformal case the Euler field assumes in the canonical coordinates
the form $E=\sum u_l \p _l$. The homogeneity relation $L_E R^{(0)}=-R^{(0)}$
together with the above flatness condition form a remarkable system of
Hamiltonian differential equations which determines $R^{(0)}$. Namely,
following B. Dubrovin \cite{Db} consider the anti-symmetric matrix 
$V_{ij}:= R^{(0)}_{ij}(u_i-u_j)$ as a point in the Poisson manifold 
$so_N^*$. Introduce $N$ non-autonomous quadratic hamiltonians 
\[ H_i:= \frac{1}{2} L_E R^{(0)}_{ii}=
\frac{1}{2}\sum_{j\neq i}\frac{V_{ij}V_{ji}}{u_i-u_j}  .\]
These hamiltonians Poisson-commute on $so_N^*$ and their flows determine
the dependence of $V$ on $u$: $\p_i V_{jl}=\{ H_i, V_{jl} \} $.
\footnote{ In fact $H_i=-R^{(0)}_{ii}/2+const$ 
in the case of conformal Frobenius 
structures. I am thankful to B. Dubrovin for this observation.
The last relation fails however in the more general setting we will
encounter in the next section.}

\medskip 

{\bf Elliptic GW-invariants.} Introduce the {\em elliptic} GW-potential of $X$
\[ G(t,q):=\sum_{d\in \L} \sum_{n=0}^{\infty} q^d [t,...,t]_n^d/n! .\]
The degree $0$ part of $G$ equals
\[ G(t,0)=-\frac{1}{24}\int_X t\w c_{\dim X -1} (T_X), \]
while the non-zero degree terms depend on $q_i$ only in the combinations
$q_i\exp t_i$ due to the divisor equation. Thus $dG$
can be considered as a closed $1$-form on the Frobenius manifold
$H=H^*(X,\CC)/2\pi i H^2(X,\ZZ )$. It has homogeneity degree $0$ with
respect to the Euler vector field on $H$.

Let $t\in H$ be a semi-simple point. In a neighborhood of $t$ the functions 
$\log \D_{\a},\ R^{(0)}_{\a\a}$ are uniquely defined up to additive constants. 
We normalize the constants in $R^{(0)}_{\a\a}$ by the homogeneity condition
\[ L_E R^{(0)}_{\a\a} = -R^{(0)}_{\a\a} .\]

\medskip

{\bf Conjecture 1.2.} 
{\em Suppose that the Frobenius structure defined by the 
rational GW-potential $F$ of the manifold $X$ is semi-simple. Then the 
elliptic GW-potential of $X$ is determined by}
\[ dG = \sum_{\a} d(\log \D_{\a})/48 + \sum_{\a} R^{(0)}_{\a\a} du_{\a}/2.\]
 
\medskip

{\em Example: $X=\CC P^1$.} The classes $1$ and $p$ Poincare-dual to the 
fundamental class and a point form a basis in $H^*(X)$. The elliptic 
GW-potential is $G=-t/24$ where $t$ is the coordinate on $H^2(X)$ since
non-constant elliptic curves do not contribute to $G$ for dimensional reasons.
Our conjecture agrees with this fact. Indeed, looking for the fundamental
solution $S=\Psi (1+R\h +o(\h))\exp (U/\h )$ of the differential system
$\h \dot{S}_1=S_2,\ \h \dot{S_2}=e^t S_1$ corresponding to the quantum 
cohomology algebra $\QQ [p,q]/(p^2-q)$ of $\CC P^1$, we find
\[ \Psi ^{-1} \dot{\Psi }=[\dot{U}, R],\ 
(\dot{R} + \Psi^{-1}\dot{\Psi} R)_{diag}=0.\]
The normalized eigen-vectors of the quantum 
multiplication operator $p\circ $ corresponding to the eigen-values 
$\dot{u}_{\pm}=\pm e^{t/2}$ are equal to $(e^{-t/4}p \pm e^{t/4})$.  
We find $(\Psi^{-1}\dot{\Psi} )_{+-}=-1/4$ and
$R_{+-}=e^{-t/2}/8=-R_{-+}$. Respectively,
\[ \dot{R}_{++}=-R_{+-}(\dot{u}_{+}-\dot{u}_{-})R_{-+}=\exp (-t/2)/32 
=-\dot{R}_{--} \]
and therefore $R_{++}du_{+}+R_{--}du_{-}=-dt/8$. 

On the other hand, $\D_{\pm}=\pm 2e^{t/2}$ and thus $d\log \D_+\D_-=dt$. 
We find $dG=dt/48-dt/16=-dt/24$. 

\section{Equivariant GW-invariants \\ in genus $0$ and $1$}

In this section we formulate a theorem confirming an equivariant version of 
Conjecture $1.2$ in the case of toric actions with isolated fixed points.

\medskip

{\bf Equivariant cohomology.}
Let a compact group $G$ act on a topological space $M$. The {\em equivariant}
cohomology $H^*_G(M)$ is defined as the cohomology of the homotopy quotient
$M_G:=(EG\times M)/G$ where $EG$ is the total space of the universal principal
$G$-bundle $EG\to BG$. When $M$ is a point $M_G=BG$, and the ring 
$H^*_G(pt)=H^*(BG)$ plays the role of the coefficient ring in the 
$G$-equivariant cohomology theory. In particular, the $G$-equvariant map
$M\to pt$ induces the $M$-bundle $M_G\to BG$ and a natural structure of the 
$H^*_G(pt)$-module on $H^*_G(M)$.

A $G$-equivariant vector bundle $V$ over a $G$-space $M$ induces a vector 
bundle $V_G$ over the homotopy quotient $M_G$. Equivariant characteristic 
classes of $V$ are defined as the ordinary characteristic classes of $V_G$.
This construction applies to equivariant orbi-bundles over orbi-spaces and
gives rise to equivariant characteristic classes of orbi-bunles well-defined
in $H^*_G(M,\QQ )$.

In the case of smooth orientation-preserving $G$-actions on compact oriented 
manifolds the {\em fiberwise
integration} over the fibres of the $M$-bundle $M_G\to BG$ defines the
$H^*_G(pt)$-linear homomorphism $\int_M: H^*_G(M)\to H^*_G(pt)$ and the
bilinear {\rm Poincar\'e pairing}
\[ \lan \phi , \psi \ran := \int_{[M]} \phi \w \psi , \]
non-degenerate over $H^*(BG,\QQ )$ in the case of Hamiltonian $G$-actions
on compact symplectic manifolds.
The same operations are well-defined over $\QQ $ in the case of orbifolds.

Let us assume now that $G$ is a torus.
According to the Borel fixed point localization formula
\[ \int_{[M]} \phi = \int_{[M^G]} \frac{i^*\phi}{Euler_G (N_M (M_G))} \]
where $i:M^G \to M$ is the inclusion of the fixed point submanifold $M^G$ into 
$M$, $N_M(M^G)$ is the normal bundle to the fixed point submanifold, 
and $Euler_G$ is the $G$-equivariant Euler class.
In the orbifold case $N_M M^G$ is an equivariant orbi-bundle over the 
orbifold $M^G$, and the localization formula holds true over $\QQ $. One should
have in mind however that the fundamental class of the fixed point 
{\em sub}-orbifold $M^G$ differs from the geometrical fundamental class
of the orbifold $M^G$ by the factors $1/|Aut |$ on each connected component,
where $Aut $ is the subgroup --- in the
symmetry group defining the orbifold structure on $M$ at a generic point of 
the component --- stabilizing the point.  

\medskip

{\bf Genus $0$.}
Let $X$ be a compact K\"ahler manifold provided with a Hamiltonian Killing
action of a compact Lie group $G$. Then the group acts also on the moduli
spaces of stable maps $X_{g,n,d}$, and this action commutes with evaluation,
forgetting and contraction maps. The constructions \cite{BF, LT, R1} of the 
virtual fundamental cycles $[X_{g,n,d}]$ can be extended to the equivariant 
setting. 
\footnote{ For convex $X$ the $G$-equivariant
virtual fundamental class $[X_{0,n,d}]$ coincides with the equivariant 
fundamental class of $X_{0,n,d}$ considered as an orbifold.}
This allows one to generalize GW-theory to the equivariant case. The theory
of equivariant genus $0$ GW-invariants \cite{Gi1} is quite analogous to the
theory of Frobenius structures reviewed in Section $1$. We describe below 
the modifications to be made in the equivariant setting emphasizing the case
of tori actions.

(1) The coefficient algebra $H^*_G(pt,\QQ)$ of the equivariant cohomology
theory replaces the ground field $\QQ $ of the non-equivariant GW-theory.
If $G$ is the $l$-dimensional torus the algebra is isomorphic to the
polynomial ring $\QQ [\l_1,...,\l_l]$ in $l$ generators of degree $1$ 
(in our complex grading) since $BG$ is weakly homotopy equivalent to
$(\CC P^{\infty})^l$. In all questions involving Borel localization formulas
the algebra is replaced by the field of fractions $\QQ (\l)$ since rational
functions of $\l $ can occur. 

(2) The equivariant GW-invariants $(T_1,...,T_n)_n^d$ and their higher genus
counterparts take values in $\QQ [\l]$ and are polylinear over $\QQ [\l ]$.
The potentials $F(t,q), S_{\a\b}(t,q,\h), V_{\a\b}(t,q,x,y)$, etc., can be 
therefore considered as formal functions with coefficients in $\QQ [\l]$ and
thus are functions of $\l$ as well. Note that $t=\sum t_{\a}\phi_{\a}$ is now
the general {\em equivariant} cohomology class of $X$, and that 
$\{ \phi_{\a} \}$ here represents a $\QQ (\l )$-basis in the 
$\QQ [\l]$-module $H^*_G(X)$.

(3) The string equation remains unchanged in the equivariant case. 
Assuming for simplicity that $X$ is simply connected (which is 
automatically the case if the Hamiltonian torus action has only isolated 
fixed points) we have the short exact sequence
\[ 0\to H^2(BG)\to H^2_G(X)\to H^2(X)\to 0 .\]
The divisor equation and its consequences for GW-potentials hold true in the 
equivariant case as well if only 
we interpret $\lan p,d\ran $ as the value of the projection of $p\in H^2_G(X)$
to $H^2(X)$ on the homology class $d\in H_2(X)$. Nevertheless it will be
convenient sometimes to keep the formal variables $q_1,...,q_r$ in place
and thus to consider $(\QQ [\l])[[\L ]]$ as the ground algebra
in the equivariant setting.  

(4) The equivariant GW-potential $F$ defines on $H:=H^*_G(X,\QQ [[\L]])$ 
a Frobenius structure over the ground ring $(\QQ [\l ])[[\L]]$.
The divisor equation induces however the following symmetry:
\[ (\p_i - q_i\p /\p q_i) \tilde{F}_{\a\b\c} (t,q,\l) =0,\ \ i=1,...,r, \]
where $\p_i$ are the tangent vector fields along $H^2_G(X)$ representing the
basis in $H^2(X)$ as in Section $1$. The grading axiom should be modified
in the equivariant case: the potential $F$ is homogeneous of degree $3-\dim X$
with respect to the Euler field
\[ E=\sum_{\a} (\deg t_{\a}) t_{\a} \p_{\a} +\sum_i c_i q_i\p /\p q_i 
+\sum_j \l_j \p /\p \l_j .\]
Due to the last summand the Euler derivation is not $\QQ [\l ]$-linear, and
thus the axioms of the conformal Frobenius manifold are not satisfied. 
We will call the Frobenius structures with such a modified grading axiom
{\em quasi-conformal}.

\medskip
 
Let us assume now that the action of the torus $G$ on $X$ has only isolated 
fixed points. The $\d$-functions of the fixed points form a basis 
$\{ \phi_{\a} \} $ in the equivariant cohomology of $X$ over the field
of fractions $\QQ (\l )$. The classical equivariant cohomology algebra of $X$
is semi-simple at generic $\l$ and therefore its quantum deformation and the
corresponding quasi-conformal Frobenius manifold $H$ is genericly semi-simple
as well. Thus Proposition $1.1$ applies and gives rise to the expansion
\footnote{We will discuss it with greater detail in Section $3$ in connection
with localization formulas.}
\[ (S_{\a\b} (t,q,\h,\l)) = \Psi (1+\h R^{(0)}+o(\h))\exp (U/\h ) \ .\]
The canonical coordinates $u_{\a}$ and the diagonal entries $R^{(0)}_{\a\a}$
of the matrix $R^{(0)}$ are defined by the Frobenius structure up to additive
``constants'' which are now elements of the ground ring $(\QQ [\l])[\L ]$.
\footnote{ Dubrovin's classification of semi-simple Frobenius structures
should be modified in the quasi-conformal case as explained in \cite{Gi1}.
In particular, the Hamiltonians $H_i$ play now the role of densities in the
Poisson-commuting Hamiltonians on the affine Lie coalgebra $\hat{so}_N^*$. 
They can be also described as the Lie derivatives 
$H_i = \sum_{\a} u_{\a}(\p /\p u_{\a}) R^{(0)}_{ii}$ but are no longer 
proportional to $R^{(0)}_{ii}$ since these functions are quasi-homogeneous with
respect to the Euler field $E$ which takes on 
$\sum u_{\a} \p/\p u_{\a} + \sum \l_j \p/\p \l_j + \sum c_i q_i\p /\p q_i$
in the canonical coordinate system.}  
Moreover, the symmetry induced by the divisor equation for the potentials 
$S_{\a\b}$ implies that $R^{(0)}$ is invariant with respect to the vector 
fields $\p_i - q_i\p /\p q_i$. This allows us to normalize the additive 
constants by the condition $R^{(0)}_{\a\a}\equiv 0 \ \mod \ (q) $.

\medskip

{\bf Genus $1$.} Introduce now the equivariant genus $1$ potential
\[ G(t,q,\l):=\sum_{d\in\L }\sum_{n=0}^{\infty} q^d[t,...,t]_n^d/n! .\]
The degree $d=0$ part of $G$ can be computed by the localization formulas:
\[ G(t,0,\l)=-\frac{1}{24}\int_X t\w c_{\dim X -1}(\calt_X) =
-\frac{1}{24}\sum_{\a} t_{\a} c_{-1}^{\a}, \]
where $c_{-1}^{\a}$ is defined to be the ratio 
$c_{\dim X-1}(\calt_X)/c_{\dim X}(\calt_X)$
of the equivariant Chern classes localized to 
the fixed point $\a \in X^G$.  

\medskip

{\bf Theorem 2.1.} {\em Suppose that the complexified action of the
torus $G_{\CC}$ on the compact K\"ahler manifold $X$ has only isolated fixed 
points and isolated $1$-dimensional orbits. Then}
\[ dG = \sum_{\a} d(\log \D_{\a})/48 -\sum_{\a} c_{-1}^{\a} du_{\a}/24 +
\sum_{\a} R^{(0)}_{\a\a}du_{\a}/2 .\]

\medskip

{\em Remarks.} (1) The differential $d$ in the theorem is taken with respect
to the coordinates $t_{\a}$ on the space $H$ over the ground ring 
$(\QQ [\l])[\L]$. Thus $q$ and $\l$ are considered as constants.

(2) Redefining the additive constants in $R^{(0)}_{\a\a}$ by
\[ R_{\a} := R^{(0)}_{\a\a}-c_{1}^{\a}/12 \]
we can reformulate the theorem in the form
\[ dG=\sum d\log \D_{\a} /48 + \sum R_{\a} du_{a}/2 \]
suggested in the introduction. In the several examples we tried
both summands in the RHS have limits as $\l$ approaches $0$ and in this
limit turn into their non-equivariant counterparts. If proven to be the 
general rule, this observation would confirm the Conjecture $1.2$ for 
toric manifolds and homogeneous K\"ahler spaces.

{\em Example.} Equivariant quantum cohomology of $\CC P^1=P(\CC^2)$ with 
respect to the circle acting by $diag (e^{i\f },e^{-i\f })$ on $\CC ^2$ is
known to be isomorphic
\footnote{We write here $q$ instead of $q\exp t$ in order to 
minimize the paperwork.}
to $\QQ [p,q,\l ]/(p^2-\l^2-q)$ with the 
equivariant Poincar\'e pairing
\[ \lan \phi , \psi \ran =\frac{1}{2\pi i} \oint 
\phi (p,q,\l ) \psi (p,q,\l ) \frac{dp}{p^2-\l^2-q} \ . \]
Introducing the indices $\pm $ for the two fixed points on $\CC P^1$
with the normal Euler classes $\pm 2\l $, we find that the normalized 
``Hessians'' $\D_{\pm}$ are equal to $\pm 2p /(\pm 2\l )$ where 
$p=(1+q/\l^2)^{1/2} =\l + O(q)$. The differentials of the canonical coordinates
are given by 
\[ du_{\pm}=\pm pd\log q = \pm \frac{2p^2 dp}{p^2-\l^2} \ .\]
In the basis $\phi_{\pm}=(\l \pm p)/2\l $ in $H^*_G(\CC P^1)$ the matrix
$\Psi $ of normalized eigen-vectors of quantum cup-product operators takes on
\[ \Psi = \frac{1}{2} \left[ 
\begin{array}{cc} z+z^{-1} & -z+z^{-1} \\ z-z^{-1} & -z-z^{-1} \end{array}
\right] \ , \]
where $z:=(1+q/\l^2)^{1/4}=1+... $. Respectively, 
\[ \Psi^{-1}\dot{\Psi} = \left[ \begin{array}{cc} 0 & -\dot{z}/z \\
-\dot{z}/z & 0 \end{array} \right]\ , \]
where the dot means $q d/dq$.

The diagonal entries of $R^{(0)}$ can be found by integration of 
$\pm (\dot{z})^2/z^2 (\dot{u}_{+}-\dot{u}_{-})=\pm (p^2-\l^2)/32p^5$.
Since $d\log q=(2p/(p^2-\l^2))dp$, we have
\[ R^{(0)}_{++}=-R^{(0)}_{--}=\int_{\l}^p \frac{x^2-\l^2}{16x^4}dx
= -\frac{1}{16p}+\frac{\l^2}{48p^3}+\frac{1}{24\l} \ .\]
The contribution to $\dot{G}$ via Theorem $2.1$ equals
\[ \frac{1}{2}(R^{(0)}_{++}\dot{u}_{+}+R^{(0)}_{--}\dot{u}_{-})=
-\frac{1}{16}+\frac{\l^2}{48p^2}+\frac{p}{24\l} \ .\]
The other two summands are
\[ -\frac{1}{24}(\frac{\dot{u}_{+}}{2\l}+\frac{\dot{u}_{-}}{-2\l})=
-\frac{p}{24} \]
(which together with the first one yields $-1/16 +\l^2/48p^2$) and
\[ 2\frac{\dot{\D}}{24\D }=\frac{\dot{p}}{24p}=\frac{p^2-\l^2}{48p^2} \ . \]
The total sum $-1/16+1/48=-1/24$ agrees with the known result $G=-(\log q)/24$.

\section{Fixed point localization \\ in genus $0$ and $1$}

The proof of Theorem $2.1$ is based on application of the Borel fixed
point localization formula to the equivariant virtual fundamental classes
$[X_{g,n,d}]$ of the moduli spaces of stable maps. A stable map
$f: (\S, \e) \to X$ represents a fixed point of the torus action in the moduli
space if its shift by the torus action can be compensated by 
automorphisms of the marked curve. Equivalently,

- $f(\S)$ is contained in the union of $0$- and $1$-dimensional orbits of 
$G_{\CC}$,

- $f(\e)$ is contained in the fixed point set $X^G$,

- if the map $f$ restricted to an irreducible component of $\S$ is 
constant then the image is a fixed point,

- if it is not constant then the component is isomorphic to $\CC P^1$,
carries no more than $2$ special points (which can be positioned at $0$
and/or $\infty $, and the map is a multiple cover $z\mapsto w=z^m$ onto
the closure of a $1$-dimensional orbit (which is also isomorphic to
$\CC P^1$ with $w=0, \infty $ to be the fixed points). 

The connected component of the fixed point set $X_{g,n,d}^G$ containing the
equivalence class $[f]$ can be described as the product of the Deligne-Mumford
spaces $\M _{g_i,k_i}$ factorized by a finite symmetry group of the 
combinatorial structure of the map $f$. Each factor $\M _{g_i,k_i}$ 
corresponds to a connected component $\S_i$ in $f^{-1}(X^G)$, $g_i$ is the
arithmetical genus of $\S_i$, and $k_i$ equals the total number of marked 
points and of special points $\overline{(\S - f^{-1}(X^G))} \cap 
f^{-1}(X^G)$ situated on $\S_i$. 

Application of the fixed point localization formula requires a description
of the virtual normal bundle to each fixed point component and reduces to 
integration over the Deligne-Mumford spaces. The idea to apply the localization
technique to the moduli spaces $X_{g,n,d}$ is due to M. Kontsevich \cite{Kn}
and was systematically exploited in \cite{Gi1, Gi3} and several other papers.
The description used in \cite{Kn, Gi1, Gi2} for localizations of the 
fundamental classes $[X_{0,n,d}]$, being obvious in the orbifold case, 
can be easily extended to the general ``virtual'' case. 
A rigorous justification of these localization
formulas was recently given in \cite{GP} on the basis of the 
algebraic-geometrical approach \cite{BF} to the virtual fundamental cycles.

\medskip

The idea of our proof of Theorem $2.1$ can be now described as follows.
Any fixed point of the torus $G$ action on the genus $1$ moduli spaces 
$X_{1,n,d}$ has the following combinatorial structure:

- either it is a tree walking along the skeleton of $1$-dimensional orbits in 
$X$, and $f^{-1}(X^G)$ has exactly one connected component $\S_0$ of 
arithmetical genus $1$,

- or it is a graph walking along the skeleton of $1$-dimensional orbits   
with exactly one cycle, and all irreducible components are rational.

The way how the fixed points of the first type contribute to the genus $1$ 
GW-potential $G$ via localization formulas can be compared to suitable genus
$0$ GW-invariants. As we shell see, the total contribution of all fixed point 
components with $\S_0$ mapped to the fixed point $\a$ in $X$ equals
$\log \D_{\a} /48 - c_{-1}^{\a}u_{\a}/24$. The proof is partly based on 
intersection theory in Deligne-Mumford spaces.

Contributions of the second type fixed points to $G$ are hard to compare with
genus $0$ GW-invariants directly because of the rotational symmetry of the 
cycle. By computing their contributions to the partial derivatives 
$\p _{\c}G$ instead, we distinguish a marked point in $(\S,\e)$ which carries
the class $\phi_{\c}$ and is situated on a branch of the graph approaching
the cycle of $1$-dimensional orbits at a fixed point $\a$. This breaks the
rotational symmetry of the cycle and allows us to compare localization formulas
for the cycles with those for the chains between the fixed points $\a$ and $\b$
with $\b=\a$. As we shell see, all such contributions add up to 
$R^{(0)}_{\a\a}\p_{\c}u_{\a}/2$. The proof is based on ``materialization''
of Dubrovin's structural theory of semi-simple Frobenius manifolds in terms
of fixed point localization.

\medskip

{\bf Intersection theory in $\M _{0,k}$ and $\M _{1,k}$.} 
     
Let us consider the GW-theory with a point taken on the role of the target 
space (so that the moduli spaces of stable maps are Deligne-Mumford spaces) 
and introduce the following GW-potentials:
\[ u(T):=\sum_{n=1}^{\infty} \frac{1}{n!} (1,T,...,T,1)_{n+2} \ ,  \]
\[ s(T,\h):=1+\sum_{n=1}^{\infty} \frac{1}{n!} (1,T,...,T,\frac{1}{\h-c})_{n+2}
\ , \]
\[ v(T,x,y):=\frac{1}{x+y}+\sum_{n=1}^{\infty} \frac{1}{n!}
(\frac{1}{x-c},T,...,T,\frac{1}{y-c})_{n+2} \ , \]
\[ \d (T):=\sum_{n=0}^{\infty} \frac{1}{n!} (1,1,1,T,...,T)_{n+3} \ , \]
\[ \mu (T):=\sum_{n=1}^{\infty} \frac{1}{n!} \o [T,...,T]_n \ , \]
\[ \nu (T):=\sum_{n=1}^{\infty} \frac{1}{n!} [T,...,T]_n \ . \]

In these formulas, $T=t_0+t_1c+t_2c^2+...$ is a series in one variable to be
replaced by the $1$-st Chern class $c^{(i)}$ of the universal cotangent line 
over $\M _{g,n}$ with the index $i$ depending on the position of the series $T$
in the correlator. The correlators $(...)_k$ and $[...]_k$ mean integration 
over $\M_{0,k}$ and $\M_{1,k}$ respectively. The correlator $\o [...]_k$
means integration
\[ \int_{\M_{1,k}} \o \w ... \]
against the $1$-st Chern class $\o $ of the {\em Hodge line bundle} ${\cal H}$
over $\M_{1,k}$. The fiber of this (orbi)-bundle over $[\S,\e]$ is the space 
$H^0(\S, K_{\S})$ of ``holomorphic differentials'' on the stable curve $\S$.
It is the pull-back of the Hodge line bundle over $\M_{1,1}$ by (any of) the
forgetting maps $\M_{1,k}\to \M_{1,1}$. Respectively, $\o \w \o =0$, the
class $\o $ on $\M_{1,1}$ coincides with $c^{(1)}$, and the orbi-structure
of $\M_{1,1}$ and ${\cal H}$ manifests in the well-known formula
\[ \int_{[\M_{1,1}]} \o = 1/24 \ .\]

The potentials $u,s,v,\d ,\mu ,\nu $ can be considered as functionals on the
space of formal series $T$. We will assume however that the coefficients
$t_0,t_1,t_2,...$ are elements of some formal series algebra $K [[\L ]]$
(in our applications $K=\QQ (\l)$), and that the whole series $T(c)$
can be rewritten as a formal $q$-series $\sum_{d\in \L } a_d(c) q^d$
with coefficients $a_d$ which are rational functions of $c$ regular at $c=0$.
Thus each $t_i$ is a formal $q$-series, and we will assume also, that 
$t_i \equiv 0 \ \mod (q) $ for $i>0$. These conditions (satisfied in our 
applications) guarantee that the tragectory of the vector field 
\[ {\cal L}:= \p /\p t_0 - t_1\p /\p t_0 - t_2 \p/\p t_1 - ... \]
with the initial condition $T$ is well-defined by   
\[ t_0(\t )=\t +\sum_{n=0}^{\infty}t_n(0)\frac{(-\t)^n}{n!},\ 
t_1(\t )=1-dt_0/d\t ,\ t_2(\t )=-dt_1/d\t , ... \]
and has a unique intersection with the hyperplane $t_0=0$. We will use these
facts in the following application of the string equation 
(notice that ${\cal L} T = 1 - (T(c)-T(0))/c $). 

\medskip

{\bf Proposition 3.1} (see \cite{Db, Gi1, D}). 
\[ s=e^{u/\h},\ v=\frac{e^{u/x +u/y}}{x+y},\ \mu=\frac{u}{24},\ 
\nu=\frac{\log \d}{24} \ .\]  

\medskip

{\em Proof.} The string equation implies 
\[ {\cal L}u=1,\ {\cal L}s=\frac{s}{\h},\ {\cal L}v=\frac{v}{x}+\frac{v}{y},\]
\[ {\cal L} \d =0,\ {\cal L} \mu = \frac{1}{24},\ {\cal L}\nu =0 .\]
At $t_0=0$ the initial conditions  
\[ u=0,\ s=1,\ v=1/(x+y),\ \d =1/(1-t_1), \mu=0,\ \nu=-(\log (1-t_1))/24 \]
can be computed from definitions with the use of dimensional reasoning
($\dim M_{0,n+2}=n-1<n$ and $\dim M_{1,n}=n$) and in the case of $\d $ and
$\nu $ --- on the basis of the formulas $(1,1,1,c,...,c)_{n+3}=n!(1,1,1)=n!$
and $[c,...,c]_n=(n-1)![c]=(n-1)!/24$ which follow from the famous 
{\em dilation equation} 
$\lan T,...,T,c \ran_{g,n+1,d}=(2g-2+n) \lan T,...,T \ran_{g,n,d} $.

\medskip

{\bf Materialization of canonical coordinates.} 
We have to review here some important
structural results from \cite{Gi1} on localization in genus $0$ equivariant 
GW-theory which were perhaps overshadowed by mirror theorems proved there. 
These results begin with the observation that a stable map $f: (\S,\e )\to X$
representing a fixed point of the torus
action on $X_{0,n,d}$ and carrying $k>3$ marked points {\em in a specific
generic configuration} (for example, $4$ marked points with a given generic
cross-ratio) must contain a connected component $\S_0$ of $f^{-1}(X^G)$ 
(we call it {\em special}) with
$k$ {\em special} ($=$ marked or singular) points realizing the given 
configuration. This follows from
the definition \cite{Kn} of the contraction maps $X_{0,n,d}\to \M_{0,k}$.
The special component is mapped to one of the fixed points $\a\in X$. This
allows us to partition certain GW-invariants into contributions --- 
via fixed point localization formulas --- of those fixed
points which map the special component to a given fixed point $\a$.
Applying the WDVV-argument to such contributions separately for each $\a$
we arrive at
some {\em local} WDVV-identities which are essentially independent on the
global WDVV-equation.  Combining local and global WDVV-equations we 
obtain simultaneous diagonalization of quantum cup-product operators
in a basis associated with fixed points of $G$ in $X$.

Let us introduce the local GW-potentials $u_{\a}, D_{\a}, \Psi_{\b}^{\a}$ 
involved into the diagonalizing structure. Let $t=\sum_{\a} t_{\a}\phi_{\a}$
denote the general equivariant cohomology class of $X$ represented in the
basis of fixed points. We denote $\p_{\a} $ 
the partial derivative with respect to
$t_{\a}$ and use the notation $\p_0$ for the differentiation operator
$\sum_{\a} \p_{\a}$ in the direction of $1\in H^*_G(X)$.
We put $e_{\a}:=\lan \phi_{\a},\phi_{\a} \ran^{-1}=\n^{\a\a}=
Euler_G (T_{\a}X)$.

\begin{itemize}
\item  Consider a point in $X_{0,n,d}$ with the property that the first two
marked points are located on the same connected component of 
$f^{-1}(X^G)$.
The total contribution of all such fixed points to the GW-potential 
\[ e_{\a} \p_{\a} \p_{\a} F_0 = 
e_{\a} \sum_{n=0}^{\infty} \frac{1}{n!} \sum_{d\in \L} q^d 
( \phi_{\a},\phi_{\a}, t,...,t )_{n+2}^d \]   
is denoted $u_{\a}$. The potentials $u_{\a}$ have homogeneity degree $1$,
are congruent to $t_{\a}$ modulo $(q)$ and can be taken on the role of
local coordinates on $H^*_G(X)$ instead of $t_{\a}$. 
\item Similarly, consider those fixed points where the first $3$ marked
points are located on the same connected component of $f^{-1}(X^G)$.
The total contribution of such fixed points to the GW-potential
\[ e_{\a} \p_{\a}\p_{\a}\p_{\a} F_0 =
e_{\a} \sum_{n=0}^{\infty} \frac{1}{n!}\sum_{d\in \L}
(\phi_{\a},\phi_{\a},\phi_{\a},t,...,t)_{n+3}^d \]
is denoted $D_{\a}$. We have: $\deg D_{\a}=0$ and $D_{\a}\equiv 1 \mod (q)$.
\item Finally, consider the fixed points with the $1$-st marked point
situated on the same connected component of $f^{-1}(X^G)$ as the two special
points which give birth to the branches carrying the $2$-nd and $3$-rd marked
points. The total contribution of such fixed points to
\[ e_{\a}\p_{\a}\p_{\b}\p_0 F_0 = 
e_{\a} \sum_{n=0}^{\infty}\frac{1}{n!}\sum_{d\in\L}q^d
( \phi_{\a},\phi_{\b}, 1, t,...,t)_{n+3}^d \]   
is denoted $\Psi_{\b}^{\a}$, has degree $0$ and reduces to $\d_{\b\a}$ modulo 
$(q)$.
\end{itemize}

\medskip

{\bf Theorem 3.2.} (see \cite{Gi1}). {\em   
The matrix $(\Psi_{\a}^i)$ satisfies the orthogonality
relations
\[ \sum_i \Psi_{\a}^i e_i^{-1}\Psi_{\b}^i = \d_{\a\b} e^{-1}_{\b}\ ,\  
\sum_{\a} \Psi_{\a}^i e_{\a}\Psi_{\a}^j=\d_{ij} e_j ,\]
the normalization condition
\[ \sum_{\a}\Psi_{\a}^i =D_i^{-1} \]
and diagonalizes structural constants of the quantum multiplication:
\[ \p_{\a}\p_{\b}\p_{\c} F_0 = 
\sum_{i\in X^G} \Psi_{\a}^i \frac{D_i\Psi_{\b}^i}{e_i} \Psi_{\c}^i  .\]  
The eigen-values $D_i\Psi_{\b}^i$ of the quantum cup-product operators 
$\phi_{\b} \circ$ satisfy the
integrability condition} 
\[ \sum_{\b} D_i \Psi_{\b}^i dt_{\b} = d u_i .\]

The theorem means that the local GW-potentials $u_i$ are the canonical 
coordinates of Dubrovin's axiomatic theory of Frobenius structures,
and $D_i$ are the square roots of the normalized ``Hessians'': 
$\D_i=e_i D_i^2$.

\medskip  

{\bf Proof of Theorem 2.1.} Let us begin with a remark on the general structure
of fixed point localization formulas in $X_{g,n,d}$. A connected component of
the fixed point set $X_{g,n,d}^G$ is identified by the combinatorial structure
of a stable map $f:(\S,\e)\to X$. The combinatorial structure can be specified
by the following data:

- the genera $g_i$ of connected components $\S_i$ of $f^{-1}(X^G)$ 
(we will call $\S_i$ {\em vertices} and consider the genus $g_i$ undefined
in the case if the vertex $\S_i$ is a point),

- the fixed points $\a = f(\S_i)$,

- the graph of rational components of $\S$ connecting $\S_i$'s
(we will call such components {\em edges}),

- the $1$-dimensional orbits in $X$ to which the edges
are mapped to and the multiplicities of the maps (when necessary we will
specify the orbit by the indices $\a \neq \b$ of the fixed points it connects,
denote $d_{\a\b}\in \L$ the degree of the orbit as a curve in $X$ and denote
$m$ the multiplicity of the map),

- the indices of marked points situated on each $\S_i$.

The connected components of $X_{g,n,d}^G$ are orbifolds (quotients of products
of Deligne-Mumford spaces), and the virtual normal bundles whose equivariant 
Euler classes occur in the localization formulas are orbi-bundles over these
orbifolds. These bundles can be split into virtual sums of contributions 
corresponding to the edges and to the vertices, and the Euler classes ---
into products of corresponding contributions. The contribution of each edge
to the Euler class has the form of the product of characters of $Lie G$ and
depends only on the correspondind degree $d_{\a\b}$ and multiplicity $m$.

Let us consider the intersection point $x$ of a vertex $\S_i$ with an edge 
$\CC P^1$. The virtual normal space contains the summand 
$T_x\S_i\otimes T_x\CC P^1$. 
It contributes to the {\em inverse} Euler class by 
$(\x_{\a\b}/m-c)^{-1}$. 
Here $\x_{\a\b}$ is the character of the torus action on
the tangent line to the closure of the $1$-dimensional orbit at $\a=f(\S_i)$,
and $c$ is the $1$-st Chern class of the universal cotangent line over 
$\M_{g_i,k_i}$ at the marked point corresponding to $x$. The product of such
contributions is to be integrated over $\M_{g_i,k_i}$ in the localization
formulas.

Adding up the contributions of all fixed point components in all the moduli
spaces $X_{g,n,d}$ with various $n$ and $d$ to certain local GW-potentials
we will obtain the exponential-like sums $\sum_k \lan T,...,T \ran /k!$ of 
integrals over the $\M_{g,k}$ with rather complicated (and unspecified) series
$T=t_0+t_1c+t_2c^2+...$. For example, the local genus $0$ GW-potential 
$u_{\a}$ equals $u(T)$. Here the $q$-series $T=\sum_d a_d(c) q^d$ has some
rational functions of $c$ on the role of the coefficients $a_d$ and 
$a_0=t_{\a}$. The whole series $T$ takes in account contributions via 
localization formulas of all tree-like branches of genus $0$ maps 
$f:(\S,\e)\to X$ which join the special component $\S_0$ at the
fixed point $\a$.    

\medskip

With the above remarks in mind, let us study now the contributions to
the genus $1$ GW-potential $G$ of all the first type fixed points whose
elliptic vertex $\S_0$ is mapped to $\a \in X$. 
The contribution of the vertex to the virtual normal bundle
contains the summand 
\[ H^0(\S_0, {\cal O}_{\S_0} \otimes T_{\a}X) \ominus 
H^1(\S_0, {\cal O}_{\S_0}\otimes T_{\a}X) .\]
Thus the inverse Euler class contains the factor
\[ \frac{Euler_G ({\cal H}^*\otimes T_{\a}X)}{Euler_G (T_{\a}X)} =
1 -  c_{-1}^{\a} \o \]
in addition to the factors $(\x_{\a\b}/m-c)^{-1}$ discussed above.
The total contribution of the first type fixed points equals therefore
\[ \sum [T,...,T]/k!  - c_{-1}^{\a} \sum \o [T,...,T]/k! 
= \nu (T)-c_{-1}^{\a} \mu (T) \ . \]
Notice that the series $T$ here is the same as in the above description of 
$u_{\a}$. We conclude that the total contribution equals
\[ (\log D_{\a})/24 - c_{-1}^{\a} u_{\a}/24 \ . \]

\medskip

Consider now the contributions of the second type fixed points to 
\[ \p_{\c} G = \sum_{n,d} q^d [\phi_{\c}, t,...,t]_{n+1}^d / n! \ . \]
Let $\a $ be the fixed point in $X$ where the tree-like branch of $(\S,\e)$
carrying the $1$-st marked point (with the class $\phi_{\c}$) joins the
cycle of edges in $\S$, and let $\S_0$ be the corresponding vertex of $\S$.
Denote $\x$ and $\x'$ the characters of $Lie G$ on the tangent lines to
the $1$-dimensional orbits where the edges of the cycle adjecent to $\S_0$
are mapped to, and denote $m$ and $m'$ the corresponding multiplicities.
Summing over all the second type fixed points with these data we see
that the contribution of the vertex $\S_0$ can be described as
$ e_{\a}^{-1} \p_{\c} v (T, \x/m, \x'/m')$
\[ = \p_{\c} \frac{\exp (u_{\a}m/\x+u_{\a}m'\x')}
{(\x/m+\x'/m')e_{\a}} =\frac{\exp (u_{\a}m/\x +u_{\a}m'/\x')}
{(\x/m)(\x'/m')e_{\a}}\ \p_{\c}u_{\a} .\]
This localization factor can be rewritten as
\[ (\p_{\c} u_{\a})\ \lim_{x,y\to 0} \frac{\exp (u_{\a}m/\x)}{e_{\a}(x+\x/m)}
\ e_{\a} \ \frac{\exp (u_{\a}m'/\x')}{e_{\a}(y+\x'/m')} \ . \]

Let us compare now this localization factor and the contribution of the rest
of the cycle with localization formulas for the genus $0$ GW-potential
\[ V_{\a\a}(x,y) = \sum_{n,d} \frac{q^d}{n!} (\frac{\phi_{\a}}{x-c},t,...,t,
\frac{\phi_{\a}}{y-c})_{n+2}^d \ . \]
The localization factors corresponding to vertices carrying 
the first and the last marked points (with the classes $\phi_{\a}$) 
vanish unless these vertices are mapped to $\a \in X$. If they are, consider
the chain of edges connecting the vertices and denote $m,m'$ and $\x,\x'$ the
multiplicities and the characters of the edges adjecent to these vertices. 
The localization factors of the vertices with these data are equal to
\[ e_{\a}^{-1} v(T,x,\x/m)=\frac{\exp (u_{\a}/x+u_{\a}m/\x)}{e_{\a}(x+\x/m)},\]
\[ e_{\a}^{-1} v(T,y,\x'/m')=\frac{\exp (u_{\a}/y+u_{\a}m'/\x')}
{e_{\a}(y+\x'/m')} . \]
Since the rest of the chain contributes to $\p_{\c} G$ and to $V_{\a\a}$ 
in the same way, we conclude that the total contribution to $\p_{\c} G$
of the second type fixed point in question equals
\[ (\p_{c} u_{\a}) \frac{1}{2} \lim_{x,y\to 0} \ [ 
e^{-u_{\a}/x} V_{\a\a} (x,y) e^{-u_{\a}/y} e_{\a} - \frac{1}{x+y} ] \ .\]      
where the factor $1/2$ takes care of the two orientations of cycles.

Let us look now at the fundamental solution matrix 
\[ S_{\b\a}=\sum_{n,d}\frac{q^d}{n!} 
(\phi_{\b},t,...,t,\frac{\phi_{\a}}{\h-c})_{n+2}^d \]
via localization formulas. The dependence of $S_{\b\a}$ on $\h$ is due only 
to the localization factor of the vertex carrying the last marked point; it is
equal to
\[ \frac{\exp (u_{\a}/\h )}{e_{\a} (\h +\x /m)} \]
if the first marked point belongs to another vertex, and to
$\exp (u_{\a}/\h) e_{\a}^{-1} \d_{\b\a}$ if the vertex is the same.
Since $\x \neq 0$, we can expand $(\h +\x /m)^{-1}$ into a power series in $\h$
and summing over all fixed point conponents obtain the asymptotical expansion
$ \Psi (1+\h R^{(0)}+o(\h))\exp (U/\h)$
of Proposition $1.1$ for the fundamental solution matrix $(S_{\b\a}e_{\a})$.
Notice that the matrix $\Psi $ of eigen-vectors here is normalized in the 
same way as (and thus coincides with) the matrix $(\Psi_{\b}^{\a})$ in 
Theorem $3.2$, since $S_{\b\a}e_{\a}\equiv \d_{\a\b} \ \mod \ (q) $.

It remains only to invoke the WDVV-identity 
\[ V_{\a\a}=\sum_{\b} S_{\b\a}(x) e_{\b} S_{\b\a}(y)/(x+y)\ , \]
the asymptotical expansion 
\[ S_{\b\a}(\h)=\sum_i \Psi_{\b}^{i} (\d_{i\a} +\h R^{(0)}_{i\a} + o(\h) )
\ e^{u_{\a}/\h} e_{\a}^{-1} \ , \]
and the orthogonality relation 
\[ \sum_{\b} \Psi_{\b}^i e_{\b} \Psi_{\b}^j=\d_{ij} e_i \]
in order to identify the above limit with $R^{(0)}_{\a\a}$.

\medskip

Combining the contributions of all first and second type fixed points we
conclude that
\[ dG = \sum_{\a} [ \ d(\log D_{\a})/24 - c_{-1}^{\a}du_{\a}/24 + 
R^{(0)}_{\a\a} du_{\a}/2 \ ] \ . \]

\section{A mirror theory for concave bundles.}

{\bf Genus $1$.}
Let $X$ be a compact K\"ahler manifold and $V$ be a holomorphic vector bundle
$E\to X$ with the total space $E$. We call the bundle $V$ {\em concave} if
for any non-constant stable map $f: (\S ,\e)\to X$ the induced bundle
$f^*V$ over $\S$ has no global holomorphic sections: $H^0(\S, f^*V)=0$.
Direct sums of negative line bundles are concave and will play the role of
main examples in this section.

If $V$ is concave then non-constant stable maps to $E$ are actually maps to
the zero section of $V$ and therefore the moduli spaces $E_{g,n,d}=X_{g,n,d}$ 
are compact for $d\neq 0$. 
This allows one to define GW-invariants of non-compact space $E$.
Namely, for $d\neq 0$ denote $V'_{g,n,d}$ the {\em obstruction} bundle over
$X_{g,n,d}$ formed by the spaces $H^1(\S, f^*V)$. The virtual fundamental
class $[E_{g,n,d}]$ is the cap-product of $[X_{g,n,d}]$ with the Euler class
of the obstruction bundle:
\[ \int_{[E_{g,n,d}]} \Phi = \int_{[X_{g,n,d}]} \Phi \w Euler (V'_{g,n,d}) .\]
  
In order to include the concave bundle spaces into the general framework of
GW-theory one has to extend the above formula to the case $d=0$ when the
moduli spaces are non-compact. Following \cite{Gi1} we provide $V$ with
the fiberwise circle action $U_1:E$ by unitary scalar multiplication.  
The constant maps $f:(\S ,\e )\to X$ form the fixed point set 
$X\times \M_{g,n}$ of 
$U_1$-action on $E_{g,n,0}=E\times \M_{g,n}$ with the normal bundle 
$V\otimes \CC = H^0(\S ,f^*V) $. 
We introduce $U_1$-equivariant GW-invariants of $E$ for $d\neq 0$ --- by
\[ \int_{[E_{g,n,d}]} \Phi = \int_{[X_{g,n,d}]} \Phi \w Euler_{U_1} 
(V'_{g,n,d}) , \]
and for $d=0$ --- by the localization formula
\[ \int_{[E_{g,n,d}]} \Phi := \int_{[X_{g,n,d}]} \Phi \w 
\frac{Euler_{U_1} (V'_{g,n,d})}{Euler_{U_1} (V)} .\]
The GW-invariants take values in the coefficient field $\QQ (\l)$ of the 
$U_1$-equivariant theory, but the degree $d\neq 0$ invariants are defined over 
$\QQ [\l ]$ and specialize to the non-equivariant ones at $\l =0$.
The construction immediately extends to the case of GW-theory equivariant 
with respect to an additional group $G$ acting on $E\to X$.  

\medskip

As it is shown in \cite{Gi1} 
\footnote{Strictly speaking the paper deals with the case of convex base $X$, 
but the arguments easily extend to any K\"ahler base as soon as the GW-theory 
for $X$ has been worked out.}
the genus $0$ GW-invariants of concave bundle spaces $E$ have the same 
properties as equivariant GW-invariants of compact manifolds including 
WDVV, string and divisor equations. In particular
they define on the space $H=H^*(E, \QQ (\l ) [[\L ]])$ 
provided with the intersection pairing 
\[ \lan \phi , \psi \ran = \int_X \phi \w \psi \w Euler_{U_1}^{-1} (V) \]
a quasi-conformal Frobenius structure of dimension $D=\dim_{\CC} E$.
The Euler vector field on $H$ is defined by the usual rules; in particular
$\deg \l =1$, and the $1$-st Chern class of the tangent bundle is 
$c_1(\calt_E)=c_1(\calt_X) + c_1 (V)$.

\medskip

We generalize our genus $1$ theory to concave bundle spaces $E$.
Similarly to the compact case we introduce potential
\[ G(t,q,\l ):= \sum_{n,d} \frac{q^d}{n!} [ t,...,t ]_d^n \] 
encoding genus $1$ equivariant GW-invariants of $E$. The Conjecture $1.2$ 
applies. 

Let us assume now that the base $X$ is provided with a Killing 
Hamiltonian action of a torus $T$, that the action of the complexified torus 
has only isolated $0$- and $1$-dimensional orbits, and that the action can be 
lifted to the bundle $E\to X$. The Frobenius structure defined on $H$ by the
genus $0$ GW-invariants of $E$ is genericly semi-simple. In particular, the
canonical coordinates $u_{\a}$, the ``Hessians'' $\D_{\a}$ and the asymptotical
coefficients $R_{\a\a}^{(0)}$ are defined.

The reduction of the
genus $1$ potential $G$ modulo $(q)$ equals $-\sum_{\a} c_{-1}^{\a} t_{\a}/24 $
where $\sum t_{\a} \phi_{\a}$ is the coordinate representation of the general
cohomology class $t$ in the basis of $\d$-functions of the fixed points, and
$c_{-1}^{\a}$ is the ratio $c_{\dim E -1} (T_{\a}E)/c_{\dim E} (T_{\a}E)$ 
of the $T\times U_1$-equivariant Chern classes of the tangent space to 
$E$ at the fixed point $\a$. 

The following theorem --- and its proof --- is a straightforward generalization
of Theorem $2.1$ to concave bundle spaces.

\medskip

{\bf Theorem 4.1.} 
\[ dG= \sum_{\a} [\  d(\log \D_{\a})/48 - c^{\a}_{-1} du_{\a}/24 +
R_{\a\a}^{(0)}du_{\a}/2 ] .\] 

\medskip 

{\bf Genus $0$.}
We develop now a mirror theory of concave {\em toric} bundle spaces
which in principal allows one to compute their genus $0$ GW-invariants.

According to T. Delzant, a compact symplectic toric manifold $X$ with the
Picard number $r$ can be 
obtained by the symplectic reduction of a standard $\CC ^m$ by a linear
torus action $T^r:\CC ^m$ on a sutable level of the momentum map. 
\footnote{We refer to \cite{Gi3} for a detailed discussion of combinatorics,
geometry and topology of symplectic toric manifolds.}
According to F. Kirwan, the equivariant cohomology algebra 
$H^*_{T^m}(X)$ with respect to the maximal torus action
$T^m:\CC ^m$ is generated over $\QQ [\l_1,...,\l_m]=H^*(BT^m)$
by the classes $w_1,...,w_m$ of degree $2$ Poincar\'e-dual to the 
invariant cycles obtained by the reduction of the coordinate hyperplanes
in $\CC ^m$. The generators $w_j$ can be written as linear combinations
\[ w_j=\sum_{i=1}^r p_i m_{ij}-\l_j , \ j=1,...,m, \]
in terms of some classes $p_1,...,p_r$ representing a basis in $H^2(X,\ZZ )$.
We will use this basis for labeling by $d=(d_1,...,d_r)$ the degrees of
curves in $X$ and denote $\L $ the semigroup of the degrees.

Consider the concave vector bundle $V: E\to X$ which is the direct sum of
$l$ {\em negative} line bundles over $X$. Let 
\[ v_j=\l'_j-\sum_{i=1}^r p_i l_{ij}, \ j=1,...,l, \]
be the $1$-st Chern classes of the summands
equivariant with respect to the torus $G:=T^m\times T^l$ action where the
second factor acts fiberwise on $V$ by diagonal transformations.   

\medskip

Our objective is to compute the fundamental solution matrix 
$(S_{\a\b}(t,q,\h))$ for the concave bundle space $E$ at $q=1$ and 
$t\in H^0_G(E)\oplus H^2_G(E)$. According to the string and divisor equations 
it coincides with 
\[ S_{\a\b}:=\sum_{d\in \L} q^d (\phi_{\a}, 
e^{(t_0+p\log q)/\h}\frac{\phi_{\b}}{\h-c})_d ,\]
where $\{ \phi_{\a} \} $ is the basis in $H^*_G(E)$ of $\d$-functions of the 
fixed points, $q^d=\exp (\sum d_it_i)$, $p\log q=\sum p_it_i$, $t_1,...,t_r$
are coordinates on $H^2(E)$ and $t_0$ is the coordinate on $H^0(E)$.

We introduce the formal series in $q$ and $1/\h$ with vector coefficients in 
$H^*_G(E, \QQ (\l,\l'))$ by
\[ J:= 1+\frac{1}{\h}\sum_{d\neq 0} q^d \ev_* 
\frac{Euler_G (V'_{0,1,d})}{\h-c} \]
where $\ev: E_{0,1,d} \to E$ is the evaluation map. The vector-function $J$
is related to the row sum of the matrix $(S_{\a\b})$ by
\[ \sum _{\a} S_{\a\b} = \lan J, e^{(t_0+p\log q)/\h} \phi_{\b} \ran . \]
This follows from the string equation in view of $\sum \phi_{\a} =1$.

\medskip 

We introduce the hypergeometric series $I$ in $q$ and $1/\h$ with coefficients
in $H^*_G(X, \QQ(\l,\l'))$ by the following explicit formula:
\[ I:= \sum_{d\in \L} q^d \Pi_{j=1}^l 
\frac{\Pi_{k=-\infty}^{L_j(d)-1} (v_j-k\h)}{\Pi_{k=-\infty}^{-1} (v_j-k\h)}
\ \Pi_{j=1}^m \frac{\Pi_{k=-\infty}^0 (w_j+k\h)}
{\Pi_{k=-\infty}^{D_j(d)} (w_j+k\h)} \]
where $L_j(d)=\sum_i d_i l_{ij}$, $D_j(d)=\sum_i d_i m_{ij}$.

\medskip

{\bf Theorem 4.2.} {\em Suppose that the $1$-st Chern class 
$\sum_{j=1}^m w_j + \sum_{j=1}^l v_j$ of the concave
toric bundle space $E$ is non-negative. Then 
\[ e^{(t_0+p\log q)/\h} J \ \ \text{and} \ \ e^{(t_0+p\log q)/\h} I \]
coincide up to a change of variables
\[ t_0\mapsto t_0+ f_0(q)+\sum \l_j g_j(q)+\sum \l'_j h_j(q) ,   \]
\[ \log q_i\mapsto \log q_i +f_i(q),\ i=1,...,r, \]
where $f_i, g_j, h_j$ (resp. $f_0$) are $q$-series supported at $\L -0$ of
the homogeneity degree $0$ (resp. $1$).}

\medskip

{\em Remark.} The hypothesis $c_1(\calt_E)\geq 0$ guarantees that 
$\deg q^d \geq 0$ for any $d\in \L$. The series $I$ and $J$ have the 
homogeneity degree $0$. By definition $J=1+o(\h^{-1})$. The change of variables
transforming $I$ to $J$ is determined by the asymptotics
\[ I=1+ \h^{-1} [f_0+\sum p_i f_i + \sum \l_j g_j + \sum \l'_j h_j ] + 
o(\h^{-1}) .\]  
Theorem $4.2$ is quite similar to the mirror theorem $0.2$ in \cite{Gi3} for
toric convex super-manifolds and complete intersections. However in the case
of concave bundles $V$ of dimension $l>1$ the series $I$ has the asymptotics
$1+o(\h^{-1})$ (due to the factor $\Pi_{j=1}^l (v_j -0\h)$) and thus coincides
with $J$.

\medskip

{\bf Corollary 4.3.} {\em For concave toric vector bundles $E\to X$ of 
dimension $>1$ the GW-potential $J$ of the total space $E$ coincides 
with the hypergeometric series $I$, provided that $c_1(\calt_E)\geq 0$.}
  
\medskip

Actually $I=1+O(\h^{-l})$ which allows us to derive 

\medskip

{\bf Corollary 4.4.}{\em If $l>1$ then in the small equivariant quantum 
cohomology algebra of $E$ we have $v_1\circ ... \circ v_l=\f (\pm q) v_1...v_l$
where $\pm q^d = (-1)^{\sum L_j(d)} q^d$ and
\[ \f (q) =\sum_{d: \sum L_j(d)=\sum D_j(d)} 
\frac{L_1(d)!...L_l(d)!}{D_1(d)!...D_m(d)!} q^d ,\]
while any shorter quantum product of 
degree $2$ classes coincides with classical.}

\medskip

{\em Proof.} We are going to exploit the fact, that the matrix 
$S=(S_{\a\b} e_{\b}^{-1})$
satisfies $\h \p S = (p\circ ) S$, where $\p $ is the derivative $q\p /\p q$ 
in the direction of a second degree class $p$, as the base for induction. 
Suppose that for all $p$ we have $(\h \p)^k S=(p^k\circ ) S $ and thus 
$= (p^k\circ ) + O(\h^{-1}) $ with no 
terms of positive order in $\h $. Then the row sums  
$\sum_{\a}\lan \phi_{\a}, p^k\circ \phi_{\b} \ran = \lan p^k, \phi_{\b} \ran $
are constant, and we conclude that 
\[ (\h \p )^{k+1} S=\h (\p (p^k\circ )) S + (p^k\circ p)\circ S \]
has the row sums $\lan p^k\circ p, \phi_{\b} \ran + O(\h^{-1})$.
Therefore, if the row sum of $S$ has the asymptotics 
$(1+O(\h^{-l}))e^{(t_0+p\log q)/\h }$ (as in our case) we conclude by induction
that for $k<l-1$ we have $\lan p^k\circ p, \phi_{\b} \ran = 
\lan p^{k+1},\phi_{\b} \ran $ for any $\b$ and thus $p^k\circ p=p^{k+1}$.  
In the border case $k=l-1$ we find $p^k\circ p$ from the row
sum of $(\h \p)^{k+1}$ modulo $\h^{-1}$.

In the case of the row sum $I e^{(t_0+p\log q)/\h }$ we apply polarization of
the above conclusion and consecutively differentiate in the directions 
corresponding to the classes $p=v_1,...,v_l$. The resulting series equals
$v_1...v_l \f (\pm q)$ modulo $\h^{-1}$, and thus 
$v_1\circ ... \circ v_l = v_1...v_l \f (\pm q)$.

\newpage

While the present paper was in preparation, a result equivalent to 
Theorem $4.2$ in the case of concave bundles over projective spaces was 
published in \cite{LLY}. 
\footnote{ In our lecture course at UC Berkeley \cite{TEG} 
the theorem was also stated over projective spaces.}

We outline below a proof of Theorem $4.2$ which is completely parallel to 
the proof of the mirror theorem for projective and toric complete 
intersections given in \cite{Gi1} and \cite{Gi3} respectively and, as
we explain in the footnotes, is a variant of the proof given in \cite{LLY}.

\medskip

{\bf Scheme of the proof.}  
{\em Step $1$.} Fixed point localization in $E_{0,2,d}$ gives rise to a
recursion relation for $S_{\a\b}$. Namely, introduce the formal series
$J_{\a}^{\b}$ in $q$ and $1/\h$ with coefficients in $\QQ (\l,\l')$ by
\[ J_{\a}^{\b}(q,\h ):=S_{\a\b}(q,\h)e^{-(t_0+p(\b)\log q )/\h} e_{\b} , \]
where $p(\b)$ is the localization of $p$ at the fixed point $\b \in E$, and
$e_{\b}=Euler_G (T_{\b} E)$. As a $1/\h$-series, $J_{\a}^{\b}=
\d_{\a\b}+O(\h^{-1})$ by definition, but 
$\sum_{\a} J_{\a}^{\b} = 1+o(\h^{-1})$ since the row sum of the matrix 
$J_{\a\b}$ is equal to the localization $J^{\b}$ of the vector-function $J$
at the fixed point $\b \in E$.

\medskip

{\bf Proposition 4.5.} 
\[ J_{\a}^{\b}(q,\h)=\d_{\a\b} + \sum _{d\neq 0} q^d P_{\a}^{\b}\ ^{(d)}
(\h^{-1})
+\sum_{\c\neq \b} \sum_{m=1}^{\infty} J_{\a}^{\c}(q,\x_{\c\b}/m) \frac
{ q^{m d_{\c\b}} \ Coeff_{\c}^{\b}(m)}{\h-\x_{\c\b}/m } \]
{\em where $P_{\a}^{\b}\ ^{(d)}$ are polynomials in $\h^{-1}$ with coefficients
in $\QQ (\l,\l')$ and $Coeff_{\c}^{\b}(m)$ are (known) rational functions of
$(\l,\l')$.}  

\medskip

The proof of Proposition $4.5$ \cite{Gi1, Gi3} is obtained by counting 
contributions to $J_{\a}^{\b}$ of fixed points in $E_{0,2,d}$ via localization
formulas. The tree representing such a fixed point contains the vertices
$\S_{\a}$ and $\S_{\b}$ carrying the two marked points and contains a chain 
of edges connecting these vertices. In the case when $\S_{\b}$ is a point,
the last edge in the chain (it connects the fixed points $\c$ and $\b$) 
yields the localization factor $Coeff_{\c}^{\b}(m)/ (\h-\x_{\c\b})$. 
The characteristics $m, \x_{\c\b}, d_{\c\b}$ of the edge here are the
same as in Section $3$. The rest
of the chain is taken care of by the localization factor 
$J_{\a}^{\c}(q,\x_{\c\b}/m)$. All other fixed points (with $\S_{\b}$ being a 
curve) contribute somehow to the polynomial tail $\sum P^{(d)} q^d$.  

Proposition $4.5$ means that the coefficients of the $q$-series $J_{\a}^{\b}$ 
are rational functions in $\h$ with $1$-st order pole at $\h=\x_{\c\b}/m$ and
the residue at the pole, controlled recursively by the (known) coefficients 
$Coeff_{\c}^{\b}(m)$, and with high order pole at $\h=0$.
 
The same recursion relation (with the index $\a$ omitted) holds true for the
row sums $J^{\b}$. There is no need here to write down explicitly the 
recursion coefficients $Coeff_{\c}^{\b}(m)$. In fact the localization 
components $I^{\b}$ of the explicitly written vector-function $I$ are also
$q$-series with coefficients rational in $\h$. Rewriting $I^{\b}$ as sums of
simple fractions $1/(\h-\x_{\c\b}/m)$ yields a recursion relation for $I^{\b}$
of the same form as the one for $J^{\b}$. It suffices to tell only that the
recursion coefficients $Coeff_{\c}^{\b}(m)$ for $I^{\b}$ are {\em the same} as
for $J^{\b}$ (while the polynomial tails can be different).  
 
\medskip

{\em Step $2$.} Given the polynomials $P^{(d)}_{\a\b}(\h^{-1})$, the recursion
relation of Proposition $4.5$ determines the matrix $(J_{\a}^{\b})$ 
unambiguously. The following proposition provides a serious constraint on these
polynomials.

\medskip

{\bf Proposition 4.6.}{\em For any $\a,\c$ the series in $q$ and 
$z=(z_1,...,z_r)$
\[ \sum_{\b} J_{\a\b} (qe^{\h z},\h) e^{p_{\b}z} e_{\b}^{-1} 
J_{\b\c} (q,-\h) \]
has coefficients polynomial in $\h$.}

\medskip

The proof of Proposition $4.6$ (see \cite{Gi1,Gi3}) is based on another
interpretation of GW-potentials $S_{\a\b}$. Let us consider the concave
bundle ${\cal E}\to {\cal X}$ which is the cartesian product of $E\to X$
with $\CC P^1$ and is provided with the standard action of $S^1$ via the
second factor. The $G\times S^1$-equivariant cohomology of ${\cal E}$
is isomorphic to the tensor product of $S^1$-equivariant cohomology algebra
$\QQ [\pi , \h] /(\pi (\pi -\h))$ of $\CC P^1$ (here $\h$ is the generator
of $H^*(BS^1)$) with the $G$-equivariant cohomology algebra of $E$. It has
a basis $\{ \phi_{\a} \pi , \phi_{\c} (\h -\pi ) \} $. 

Let ${\cal E}_{n,d}$ be the moduli space of genus $0$ stable
maps to ${\cal E}=E\times \CC P^1$ with $n$ marked points of degre $d$ in
projection to $E$ and of degree $1$ in projection to $\CC P^1$. We introduce 
the GW-potential
\[ {\cal G_{\a\c}}:= \sum_{n,d} \frac{q^d}{n!} 
\lan \phi_{\a} \pi, pz\pi, ... , pz\pi , \phi_{\c} (\h-\pi) \ran_d^{n+2} ,\]
where the correlator $\lan ... \ran_d^{n+2}$ refers to the 
$G\times S^1$-equivariant GW-invariant of the concave bundle space ${\cal E}$
obtained by integration over $[ {\cal E}_{n+2,d} ]$.  

The series ${\cal G}_{\a\c}$ depends on $\h$ but not on $\h^{-1}$ since it
is defined without localization to fixed points of $S^1$-action. Applying
localization to fixed points of $S^1$-action on ${\cal E}_{n,d}$ and then
using the divisor equation one finds that ${\cal G}_{\a\c}$ coincides with
the convolution series introduced in Proposition $4.6$.

{\em Remark.} One can also define the GW-potentials ${\cal G}_{\a\c}$ by
\[ {\cal G}_{\a\c}=\sum_d q^d \int_{[{\cal E}_{2,d}]}
\ev_1^*(\phi_{\a} \pi) \ev_2^*(\phi_{\c} (\h-\pi )) 
Euler_{G\times S^1} ({\cal V}'_{2,d}) e^{Pz} ,\]
where $P=(P_1,...,P_r)$ are the equivariant $1$-st Chern classes of the 
{\em universal} line bundles over ${\cal E}_{2,d}$ introduced in \cite{Gi3}.
The definition uses embeddings of the toric manifold into
projective spaces and the map $\f $ from ${\cal X}_{0,d}$
to the {\em toric compactification} ${\cal X}_{(d)}$ of spaces of degree $d$ 
maps $\CC P^1\to X$ in the case when $X$ is a projective space. 

Summing ${\cal G}_{\a\c}$ over $\a$ and $\c$ we arrive to the polynomiality
property for the row sums $J^{\b}$. The same polynomiality property holds true
for $I^{\b}$: $ \lan I(qe^{\h z},\h) , e^{pz} I(q,-\h) \ran $ depends on $\h $ 
but not on $\h^{-1}$. The proof (see \cite{Gi3}) is based on localization 
to fixed point of $S^1$-action applied to the series
\[ \tilde{\cal G}=\sum_{[{\cal X}_{(d)}]} Euler_{G\times S^1} ({\cal V'}_{(d)})
e^{Pz} \]
which mimics the GW-potential $\sum_{\a\c}{\cal G}_{\a\c}$ in terms of toric
compactifications ${\cal X}_{(d)}$ of spaces of degree $d$ maps $\CC P^1\to X$.

\medskip

{\em Step $3$.} Let us call a solution $(J_{\a}^{\b})$ 
(respectively $(J^{\b})$) to the recursion relation of Proposition $4.5$ 
{\em polynomial} if it satisfies the polynomiality property described in
Proposition $4.6$. 

\medskip

{\bf Proposition 4.7.} {\em A polynomial solution $(J_{\a}^{\b})$ 
(respectively $(J^{\b})$) to the 
recursion relation satisfying the asymptotical conditions 
\[ J_{\a}^{\b}=\d_{\a\b}+O(\h^{-1}), \ \sum_{\a} J_{\a}^{\b}=1+o(\h^{-1}) \]
(respectively $J^{\b}=1+o(\h^{-1})$) is unique (if it exisits).}

\medskip

The proof is obtained by a straightforward argument of perturbation theory
as in Proposition $4.5$ in \cite{Gi3}.
This result completes the proof of Corollary $4.3$. 

{\em Remark.} It is not hard to prove
that under the hypotheses of Corollary $4.3$ there exist differential operators
${\cal D}_{\a}(\h q\p/\p q, q, \h)$ such that 
\[ J_{\a}^{\b} e^{(p_{\b}\log q)/\h} =  {\cal D}_{\a}
[J^{\b} e^{(p_{\b}\log q)/\h}] .\]
The proof is constructive, but we do not know how to describe the formulas
for $J_{\a}^{\b}$ in a closed form.

\medskip

{\em Step $4$.} Since both $(J^{\b})$ and $(I^{\b})$ are polynomial
solutions to the same recursion relation, the proof of Theorem $4.2$ is 
completed by the following proposition whose proof is also straightforward.
  
\medskip

{\bf Proposition 4.8.} {\em Transformations described in Theorem $4.2$
preserve the class of polynomial solutions to the recursion relation.}

\medskip

{\em Remark.} One can easily point out constant coefficient differential 
operators ${\cal D}_{\a}$ such that the hypergeometric series 
 \[ I_{\a}^{\b}:= e^{(-p\log q)/\h} {\cal D}_{\a} [ I^{\b} e^{(p\log q)/\h} ]
\]
form a polynomial solution to the recursion relation of Proposition $4.5$.
They usually have wrong asymptotics however. The matrix $(I_{\a}^{\b})$ can
be transformed to $J_{\a}^{\b}$, but in general the transformation requires
matrix differential operators of infinite order 
(including changes of variables), and we do not know how to describe the 
transformation concisely.  
\footnote{The genus $0$ mirror conjecture for complete intersections 
in the projective 
space $X=\CC P^n$ has now five proofs --- 
the four variations of the same proof  
(in \cite{Gi1}, in \cite{Gi3}, the one outlined above but applied to convex 
bundles over $X$ instead of concave bundles, and the one in 
Section $5$ of this paper based on nonlinear Serre duality), 
and the proof recently given in \cite{LLY}. Here we compare the methods 
in \cite{LLY} with our approach.

The key idea (see Step $2$ above) --- 
to study GW-invariants of the product $X\times \CC P^1$ 
{\em equivariant with respect to the $S^1$-action on $\CC P^1$} instead of
GW-invariants on $X$ --- is borrowed in \cite{LLY} from our paper \cite{Gi1},
Sections $6$ and $11$. In fact this idea is profoundly rooted in the heuristic 
interpretation \cite{Gi4} of GW-invariants of $X$ in terms of Floer 
cohomology theory on the loop space $LX$ where the $S^1$-action is given 
by rotation of loops. The generator in the cohomology algebra of $BS^1$ denoted
$\hbar $ in our papers corresponds to $\a $ in \cite{LLY}.

Another idea, which is used in all known proofs and is due to M. Kontsevich
\cite{Kn}, is to replace the virtual fundamental cycles of spaces
of curves in a complete intersection by the Euler cycles of suitable vector 
bundles over spaces of curves in the ambient space. Both papers \cite{Gi1}
and \cite{LLY} are based on computing the push forward of such cycles to 
simpler spaces. Namely, the cycles are $S^1$-equivariant Euler classes of
suitable bundles over stable map compactifications of spaces of bi-degree 
$(d,1)$ rational curves in $X\times \CC P^1$, the simpler spaces are {\em toric
compactifications} of spaces of degree $d$ maps $\CC P^1 \to X=\CC P^n$,
and the push-forwards are denoted 
$E_d$ in \cite{Gi1} and $\varphi_{!}(\chi_d )$ in \cite{LLY}.

The toric compactification is just the projective space $\CC P^{(n+1)d+n}$
of $(n+1)$-tuples of degree $\leq d$ polynomials in one variable $z$, which 
genericly describe degree $d$ maps $\CC P^1\to X=\CC P^n$; the space is 
provided with the $S^1$-action $z\mapsto z \exp(it)$ (as in the loop space!) 
Thus both papers
depend on continuity of certain natural map (denoted $\mu $ in \cite{Gi1} and 
$\varphi $ in \cite{LLY}) between the two compactifications. The continuity
is stated in \cite{LLY} as Lemma $2.6$. It coincides with our Main Lemma
in \cite{Gi1}, Section $11$. The proof of Lemma $2.6$ attributed in \cite{LLY}
to J. Li coincides with our proof of the Main Lemma. 
The difference occurs in the
proof of a key step formulated as {\em Claim} in \cite{Gi1}: our proof of the
Claim by bare hand inductive computation in the spirit of G. Segal's 
representation of vector bundles over curves via loop groups is replaced
(and this is the contribution of J. Li) by a more standard 
algebraic-geometrical argument based on the proof of Theorem $9.9$ in 
Hartshorne's book. It is worth repeating here the remark from \cite{Gi1} 
that a different proof of the lemma was provided to me
by M. Kontsevich, with whom we first discussed the map between the two 
compactifications in Fall $1994$.  

The new concept introduced in \cite{LLY} --- the {\em eulerity} property
of the classes $E_d$ (Definition $2.3$ in \cite{LLY}) --- is
to replace both the {\em recursion} relation (Step $1$ above) and the 
{\em polynomiality} property (Step $2$) of the gravitational GW-invariant 
($J$ in the above outline). Eulerity is actually equivalent to recursion $+$
polynomiality. Theorem $2.5$ in \cite{LLY} asserting the eulerity property of 
the classes $\{ E_d \}$ coincides with Proposition $11.4 (2)$ in \cite{Gi1}
deduced there {\em from} the recursion $+$ polynomiality.
The proof of Theorem $2.5$ in \cite{LLY} is based on the same
localization to fixed points of $S^1$-action on spaces of curves as in our
proof of Corollary $6.2$ in \cite{Gi1} which guarantees the polynomiality.
The recursion is derived in \cite{Gi1} by further fixed point localization 
with respect to the torus acting on $X=\CC P^n$. 
Thus the proof in \cite{LLY} shows that the latter
localization argument is unnecessary. 
      
The relationship among the two solutions to the recursion relation --- the 
gravitational GW-invariant and the explicitly defined hypergeometric series 
($I$ in the above outline) --- is based on some
uniqueness result (Proposition $11.5$ in \cite{Gi1}) for solutions to the
recursion relation satisfying the polynomiality property. The corresponding
result in \cite{LLY} is Theorem $2.11$ about linked Euler data. {\em Linked}
there translates to our terminology as the recursion coefficients in the 
recursion relations for $I$ and $J$ being the same. The proof of the
uniqueness result in \cite{LLY} is the same as in \cite{Gi1} or \cite{Gi3}.
The difference is that the uniqueness property is formulated in \cite{LLY} 
solely in terms of the Euler data $\{ E_d \}$ and not in terms of 
gravitational GW-invariant the data generate. 

The uniqueness result allows to identify the gravitational and hypergeometric 
solutions to the recursion by some changes of variables (the mirror 
transformations). This is deduced in \cite{Gi1} from Proposition $11.6$ which
states that both the recursion relation and the polynomiality property are 
preserved by the mirror transformation (see Step $4$ above). 
The corresponding result in \cite{LLY} 
is Lemma $2.15$ which says that the (equivalent!) eulerity property is 
invariant under mirror transformations. It turns out however that while it is 
straightforward to check the invariance of recursion and polynomiality 
(Proposition $11.6$ in \cite{Gi1}), it is technically harder to give 
a direct proof of the invariance of eulerity, which requires the notion of 
{\em lagrangian lifts} introduced in \cite{LLY}. The use of lagrangian lifts
is therefore unnecessary.

The last part of the proof in \cite{LLY} (see Section $3$ there) addresses the 
following issue: while the previous results allow to compute {\em some} 
GW-invariants in terms of hypergeometric functions, what do these 
GW-invariants have to do with the structural constants of quantum cohomology 
algebra involved in the formulation of the mirror conjecture? 

The computational approach to the issue in \cite{LLY} is also not 
free of overlaps with \cite{Gi1}. However
it remains unclear {\em to us} why the authors of \cite{LLY} ignore 
the fundamental relationship between the gravitational GW-invariant and 
quantum cohomology which resolves the issue momentarily. 
The relationship was described 
by R. Dijkgraaf and B. Dubrovin \cite{Db} in the axiomatic context of 
$2$-dimensional field theories and adjusted to
the setting of equivariant GW-theory in Section $6$ of \cite{Gi1}. According
to these results the structural constants of quantum cohomology algebra 
(such as {\em Yukawa coupling} in the case of quintic $3$-folds) are
{\em coefficients} of the linear differential equations satisfied by 
{\em the} gravitational GW-invariants in question. 
In fact such a relationship was the initial point of the whole project started 
by \cite{Gi2,Gi4} and completed in \cite{Gi1,Gi3}.

Thus the two proofs of the same theorem appear to be variants of
the same proof rather than two different ones, except that our reference
to the general theory of equivariant quantum cohomology, developed in 
\cite{Gi1}, Sections $1$ -- $6$, for concave and convex 
vector bundles over convex manifolds, is replaced in \cite{LLY} by a 
computation. 

It is worth straightening some inaccuracy of \cite{LLY} in quotation.
As it is commonly known, ``Givental's idea of studying equivariant Euler 
classes'' (see p. $1$ in \cite{LLY}) is due to M. Kontsevich \cite{Kn} 
who proposed a 
fixed point computation of such classes via summation over trees.
The idea of the equivariant version of quantum cohomology listed 
on p. $6$ of \cite{LLY} among ``a number of beautiful 
ideas introduced by Givental in \cite{Gi1, Gi2}'' was actually suggested 
two years earlier in \cite{GK} by a different group of authors.
The statement in the abstract that the paper \cite{LLY} ``is
completing the program started by Candelas et al, Kontsevich, Manin and
Givental, to compute rigorously the instanton prepotential function for
the quintic in $P^4$'' is also misleading: the paper is more likely to 
confirm that the program has been complete for two years.
}
\medskip

{\bf Mirrors.} 
The hypergeometric series $I^{\b}$ can be represented by hypergeometric
integrals:
\[ I^{\b}(q,\h) e^{(p_{\b}\log q)/\h}=\int_{\Gamma_q^{\b} \subset E'_q} 
e^{(\sum_{j=1}^m W_j + \sum_{j=1}^l V_j)/\h} 
\Pi_{j=1}^m W_j^{\l_j/\h} \Pi_{j=1}^l V_j^{-\l'_j/\h} \ \times \]
\[ \times \  \frac{d\log W_1\w ...\w d\log W_m\w d\log V_1\w ...\w d\log V_l}
{d\log q_1 \w ... \w d\log q_r} .\]
Here $\Gamma_q^{\b}$ are suitable non-compact cycles of middle dimension in
the complex $m+l-r$-dimensional manifold 
\[ E'_q = \{ (W,V) | \Pi_{j=1}^m W_j^{m_{ij}} =q_i \Pi_{j=1}^l V_j^{l_{ij}},
\ i=1,...,r \}  \]
provided with the local coefficient system $W^{\l/\h}V^{-\l'/\h}$.

Due to Theorem $4.2$ the above oscillating integral can be considered as the 
{\em mirror partner} of the concave toric bundle space $E$ 
\footnote{The series $I(q,\h)e^{(p\log q)/\h}$ is annihilated by any linear
differential operator ${\cal D} (\h q\p /\p q, q, \h)$ which annihilates the
integral with any cycle $\Gamma $, but usually not vice versa. In order to
get a one-to-one correspondence here one should impose some constraint on 
the cycles. We do not know however an exact description of the corresponding 
homology group. Different choices of such a constraint should correspond 
to different
toric bundles $E\to X$ whith the same matrices $(m_{ij})$ and $(l_{ij})$.
Thus the integral formula itself (which depends only on these matrices)
can represent mirror partners of several different spaces (depending 
on the level of the momentum map in the symplectic reduction procedure).} 
in the sense of 
the generalized mirror conjecture suggested in \cite{Gi2}: the equivariant
GW-potential $J^{\b}(q,\h) e^{(p\log q/\h)}$ which plays the role of 
oscillating integrals in our symplectic topology --- singularity theory
dictionary coincides with the oscillating integral after the transformation
to flat coordinates described in the theorem.

Regardless of the mirror theory one can use the integral representation in
order to compute the genus $1$ GW-potential of $E$ via the Hessians $\D_{\a}$ 
and the asymptotical coefficients $R_{\a}$ at the critical points of the phase
function $\sum (W_j + \l_j \log W_j) + \sum (V_j-\l'_j \log V_j)$ under the
constraints $\sum m_{ij} \log W_j - \sum l_{ij} \log V_j = \log q_i,\ 
i=1,...,r$. 
\footnote{Notice that choosing $p_1,...,p_r$ on the role of Lagrange 
multipliers we arrive at the equations of the critical points in the form 
$ W_j=\sum_i p_i m_{ij} -\l_j,\ V_j=\l'_j-\sum_i p_i l_{ij} $.}  
We suggest the reader to recover the genus $1$ potential $dG = dq/24q$ for 
$E=X=\CC P^1$ by this method and observe that the asymptotical coefficients
$R_{\a}$ coincide with $R^{(0)}_{\a\a}-c_{-1}^{\a}/12$. 

\medskip

{\bf Application.} Consider a generic holomorphic sphere $\CC P^1$ in a 
Calabi-Yau $3$-fold. Such spheres occur in a discrete fashion with the normal
bundle isomorphic to $\calo (-1)\oplus \calo (-1)$ which is concave. 
Multiple covers of 
this sphere contribute to genus $0$ and $1$ GW-potentials of the $3$-fold, and 
the problem of computing these contributions reduces to studying GW-invariants
of the non-compact total space $E$ of the normal bundle.  

According to Corollary $4.3$ the hypergeometric series
\footnote{We reduce the group $G$ here to the one-dimensional torus so that
$H^*_G(\CC P^1)=\QQ [p,\l]/(p^2-\l^2)$ where $p$ is the equivariant Chern class
of $\calo (1)$.}
\[ I=\sum_{d=0}^{\infty } q^d \frac{\Pi_{m=0}^{d-1} (p+m\h)^2}
{\Pi_{m=1}^d (p-\l+m\h)(p+\l+m\h)} \]
coincides with the GW-potential $J$. The intersection index
\[ \int_{[E]} \f (p) I = \frac{1}{2\pi i} 
\oint \f (p) I \frac{dp}{\l^2 (p^2-\l^2)} \]
equals the sum of a $(d=0)$-term which has no limit at $\l=0$ with a series
which at $\l=0$ turns into
\[ \sum_{d=1}^{\infty} q^d \Res_{p=\infty }\frac{ \f (p) dp}{p^2 (p+d\h)^2} .\]
This formula with $\f = \exp{(p\log q)/\h}$ determines the contribution of
multiple covers of the sphere to the genus $0$ GW-potentials $S_{\a\b}$ 
of the Calabi-Yau $3$-fold. The result coincides --- and Corollary $4.3$ 
explaines why --- with the formula for such a contribution obtained in 
\cite{Gi4} by toric (and therefore heuristic) methods. 
It is shown in \cite{Gi4} how this result implies the famous formula 
\footnote{The formula claimed by physicists \cite{COG} was first confirmed
in \cite{AM} by toric methods and then rigorously justified by Yu. Manin, J. 
Bryan, R. Pandharipande by equivariant methods more elementary than the mirror
theory. This example is also contained in \cite{LLY} and \cite{TEG}.}
$D^3 Q^D/(1-Q^D)$ for contributions of degree $D$ spheres and their 
multiple covers to the Yukawa 
coupling $\lan P\circ P, P\ran $ of a Calabi-Yau $3$-fold. 

Furthermore, according to Corollary $4.4$ the class $p$ satisfies the relation
$p^2=\l^2/(1-q)$. 
in the equivariant quantum cohomology algebra of the concave space $E$. 
From this, we find the differentials of the canonical coordinates: 
\[ du_{\pm}=dt_0+p_{\pm}d\log q = dt_0 \pm \l (1-q)^{-1/2} d\log q \]
and thus $ \p /\p u_+ - \p /\p u_- = (1-q)^{1/2}\l ^{-1} q\p /\p q $.
From the intersection pairing $\lan 1 , 1 \ran = \lan p ,p \ran =0$,
$\lan 1, p\ran = \l^{-2}$ in the equivariant cohomology of $E$ we find
the ``Hessians'': $1/\D_+ +1/\D_-=0, \ p_+/\D_+ + p_-/\D_-=\l^{-2}$ and thus
\[ \D_{\pm}=\pm 2\l^3 (1-q)^{-1/2} .\]
In particular $d\log (\D_+\D_-)=q(1-q)^{-1} d\log q$.
Using the expression of $R_{ij}$ in terms of $\D_i$ we find, 
after some elementary computations, 
\[ dR_{++}=-dR_{--}=\frac{q^2 d\log q}{32 \l (1-q)^{3/2}} \]
and therefore 
\[ R_{++}=-R_{--}=\frac{(1-q)^{1/2}}{16\l}+\frac{(1-q)^{-1/2}}{16\l } -
\frac{1}{8\l} .\]
Since the constants $-c_{-1}^{\pm}/12=\pm 1/8\l $, we conclude that
\[ (R_{++} - c_{-1}^+/12)du_+ + (R_{--}-c_{-1}^-/12)du_-=
(\frac{1}{8}+\frac{1}{8(1-q)}) d\log q . \]
Combining $1/2$ of this with $1/48$-th of $q(1-q)^{-1} d\log q$ we finally
arrive at
\[ dG = (\frac{1}{8}+ \frac{q}{12(1-q)}) d\log q .\]

In the application to counting multiple elliptic covers of $\CC P^1$ the 
degree $d=0$ term $1/8$ is to be ignored. The remaining part 
$d \log (1-q)^{-1/12}$  gives rise to the formula
$\log (1-Q^D)^{-1/12}$ for the contribution of degree $D$ rational curves to
the genus $1$ GW-potential of a Calabi-Yau $3$-fold, claimed by physicists 
and recently confirmed by more elementary equivariant methods in \cite{GP}. 
   
\section{Nonlinear Serre duality.}

Let $Y\subset X$ be a submanifold given by a section of a vector bundle 
$V$. It is plausible that some GW-invariants of $Y$ depend only on the 
bundle. In higher genus realization of this idea encounters some 
obstruction avoidable in the genus $0$ case which we begin with.

\medskip

{\bf Convex super-manifolds.} 
The vector bundle $V: E\to X$ is called 
{\em convex} if it is spanned by global holomorphic sections. 
For a stable genus $0$ map $f:\S \to X$ we have $H^1(\S ,f^*V)=0$ and thus
the spaces $H^0(\S, f^*V)$ form an orbi-bundle $V_{0,n,d}: E_{0,n,d}\to 
X_{0,n,d}$ over the moduli space. We introduce genus $0$ GW-invariants of the 
{\em super-manifold} $\Pi E$ by defining the virtual fundamental class
$[\Pi E_{0,n,d}] $ as the cap-product of the homology class $[X_{0,n,d}]$
with the cohomology class $Euler_G(V_{0,n,d})$. Here $Euler_G$ means the
equivariant Euler class with respect to a (lifted to $V_{0,n,d}$) Hamiltonian
action of $G$ on $E\to X$ such that all fixed points are contained in the 
zero section. Therefore the GW-invariants $(A,B,...,C)_n^d$ of the 
supermanifold $\Pi E$ take their values in $H^*(BG,\QQ)$.  

The genus $0$ equivariant GW-theory extends, as it is shown in \cite{Gi1}, 
to super-manifolds without any serious changes and gives rise to a 
Frobenius structure over the ground ring $\QQ (\l ) [[\L ]]$ on the 
equivariant cohomology space $H=H_G^*(X, \QQ [[\L]])$. The Poincar\'e metric
on $H$ is induced by the pairing
\[ \lan \phi , \psi \ran := \int _{[X]} \phi \w \psi \w Euler_G (V) .\]
The quasi-conformal structure is determined by the usual grading on $H^*_G(X)$
and by the $1$-st Chern class of the super-manifold:
\[ c_1 (\calt_{\Pi E}) := c_1(\calt_X) -c_1(V) .\]
The conformal dimension of the Frobenius structure equals 
$\dim_{\CC} X-\dim_{\CC} (V)$ which coincides with the super-dimension of 
$\Pi E$.

The holomorphic section $s: X \to E$ restricted to a curve $f:\S \to X$ induces
an element in $H^0(\S , f^*V)$ and thus --- a section 
$s_{0,n,d}: X_{0,n,d}\to E_{0,n,d}$. The zero locus $s_{0,n,d}^{-1}(0)$ 
coincides with the moduli space $Y_{0,n,d}$ of stable maps to $Y=s^{-1}(0)$
of degree $d$ {\em in the ambient space} $X$. The virtual fundamental class
$[Y_{0,n,d}]$ in $X_{0,n,d}$ coincides with $[X_{0,n,d}]\cap Euler(V_{0,n,d})$.
This (M. Kontsevich's) observation serves as a basis for
applications of GW-theory of super-manifolds to complete intersections. 
It shows that in the non-equivariant limit $\l =0$ the equivariant correlators
$(A,B,...,C)_n^d$ of $\Pi E$ turn into the corresponding correlators of $Y$ 
{\em among classes induced from the ambient space} $X$. 

\medskip

{\bf Localization via materialization.}
Let us assume now that the group $G_{\CC }$ is a torus acting on $X$ with 
isolated zero- and one-dimensional orbits, that the bundle $V$ is the sum
of positive line bundles (so that it is convex and the dual bundle $V^*$
is concave), and that the action is lifted to $V$ and $V^*$ in the dual 
fashion. We will show (following Section $12$ in \cite{Gi1}) that the
Frobenius structures of the convex super-manifold $\Pi E$ and of the concave 
bundle space $E^*$ are closely related.

Let $\phi_{\a}$, as usually, be the basis of fixed points in $H^*_G(X)$.
Denote $e^0_{\a}$ and $e'_{\a}$ the Euler factors $Euler_G (T_{\a}X)$ and
$Euler_G (V_{\a})=(-1)^{\dim V} Euler_G (V^*)$ respectively.
We put $s_{\a}^{\b}(\h) :=S_{\a\b}(\h)e^{-u_{\b}/\h} e^0_{\b} (e'_{\b})^{-1}$ 
where $u_{\b}$ are materialized canonical coordinates of Theorem $3.2$ 
applied to the case of the super-manifold $\Pi E$. 

\medskip

{\bf Proposition 5.1.} {\em The matrix $(s_{\a}^{\b})$
satisfies the recursion relation}
\[ s_{\a}^{\b} (\h) =\d_{\a\b} + \sum_{(\c ,\b )} \sum_{m=1}^{\infty} 
s_{\a}^{\c} (\x_{\c\b}/m) \ Coeff_{\c}^{\b}(m)
\frac{(q^{d_{\c\b}}e^{(u_{\c}-u_{\b})/\x_{\c\b}})^m}{\h -\x_{\c\b}/m} .\]

\medskip

The summation indices $(\c ,\b )$ indicate one-dimensional orbits of 
$G_{\CC }$ connecting the fixed points $\c$ and $\b$, and $d_{\c\b}$ and
$\x_{\c\b}$ have the same meaning as in Section $3$. 

{\em Proof.} The proposition is obtained by localization technique in the
same way as the recursion relation in Step $1$ in the proof of Theorem $4.2$:
we cut out the last edge ($m$-multiple cover of the orbit $(\c,\b)$) 
in the chain connecting the vertices carrying $\phi_{\a}$ and $\phi_{\b}$. 
The recursion coefficient takes in account the localization factor of the 
edge. However this time we use the localization factor of the last vertex 
(carrying $\phi_{\b}$) in the form described in Proposition $3.1$:  
$ v=(\h +\x_{\b\c}/m)^{-1}\exp (u_{\b}/\h+mu_{\b}/\x_{\b\c}) $. $\ \square$ 

The linear recursion relation of Proposition $5.1$ unambiguously determines 
the fundamental solution $(S_{\a\b})$ as a function of canonical coordinates.
The relation between canonical coordinates $u_{\b}$ and the flat coordinates
$t_{\b}$ is non-linear and is obtained from the asymptotics 
$\sum_{\a} S_{\a\b}=(1+o(\h^{-1}))e^{t_{\b}/\h}(e^0_{\b})^{-1} e'_{\b}$. 
Expand $s_{\a}^{\b}(\h)$ as $\d_{\a\b}+\dot{s}\ _{\a}^{\b}\h^{-1}+o(\h^{-1})$.
With this notation we arrive at
\[ t_{\b}=u_{\b}+\sum_{\a} \dot{s}_{\a}^{\b}  .\]

\medskip

Parallel results for concave bundles $V^*$ look as follows. Denote here the
fundamental solution matrix by $(S^*_{\a\b})$ and put 
$s^*\ _{\a}^{\b} = (-1)^{\dim V} e'_{a} S^*_{\a\b} e^{-u_{\b}/\h} e^0_{\b}$.
Then $(s^*\ _{\a}^{\b} (\h))$ satisfies the recursion relation of Proposition
$5.1$ with new recursion coefficients $Coeff^*\ _{\c}^{\b}(m)$, where $u_{\b}$
are now the canonical coordinates of Theorem $3.1$ applied to the concave 
bundle space $E^*$. Respectively, the flat coordinates $t_{\b}^*$ are 
found from $\sum_{\a} S^*_{\a\b}=
(1+o(\h^{-1})) e^{t_{\b}/\h} (e^0_{\b} e'_{\b})^{-1} (-1)^{\dim V}$:
\[ t^*_{\b}= u_{\b} +\sum_{\a} (e'_{\a})^{-1} \dot{s}^*\ _{\a}^{\b} e'_{\b} .\]
    
\medskip

Now the Serre duality enters the game as the following identity:
\[ Coeff^*\ _{\c}^{\b} (m) = (-1)^{m c_{\c\b}} Coeff_{\c}^{\b}(m) ,\]
where $c_{\c\b}$ is the value of the $1$-st Chern class of the bundle
$V$ on the degree $d_{\c\b}$ of the $1$-dimensional orbit $(\c\b)$.

\medskip

Indeed, the coeefficients arise from fixed point localization formulas 
applied to the map $\f : \S_0 \to \CC P^1,\ z \mapsto w=z^m $ of 
$\S_0 \simeq \CC P^1$ onto the orbit. 
Both coefficients are ratios of two equivariant Euler classes which
come from the virtual normal space to $X_{0,n,d}^G$ in $X_{0,n,d}$ (the
denominators), and from
the bundles $V_{0,n,d}$ and $V^{*'}_{0,n,d}$ respectively (the numerators). 
Due to our choice of normalization for $s$ and $s^*$ the denominators are 
the same, and the numerators are the equivariant Euler classes respectively of
the space (of holomorphic sections vanishing at $z=\infty$)
\[ H^0(\S_0, (\f ^*V) \otimes \calo_{\S_0} (-[\infty ])\ )  \]
and of
\[ H^1(\S_0, (\f ^*V^*) \otimes \calo_{\S_0} (-[0])\ ) .\]
By (elementary) Serre duality on $\S_0$ the second space is canonically dual
to
\[ H^0(\S_0, (\f ^*V) \otimes {\cal K}_{\S_0} ([0])\ ) .\]
But the twisted canonical line bundle ${\cal K}_{\S_0} ([0]+[\infty])$ on 
$\S_0 \simeq \CC P^1$  is trivialized by the invariant section $d\log z$. 
Since $mc_{\c\b}$ is the dimension of the dual cohomology spaces, 
we conclude that the numerators differ by the sign $(-1)^{mc_{\c\b}}$.

Thus we arrive at the ``nonlinear Serre duality'' theorem \cite{Gi1}.

\medskip

{\bf Theorem 5.2.} {\em The fundamental solution matrices $(S_{\a\b})$ and
$(S^*_{\a\b})$ of the dual convex super-manifold $\Pi E$ and concave vector
bundle space $E^*$, considered as functions of canonical coordinates, satisfy 
\[ S_{\a\b}(u, q,\h )=(-1)^{\dim V} e'_{\a} S^*_{\a\b} (u, \pm q, \h ) 
e'_{\b} ,\]
where $\pm q^d$ means $(-1)^{\lan c_1(V),d \ran } q^d$.}

\medskip

In flat coordinates, the fundamental solution matrices are related therefore 
by an additional transformation of coordinates $t^*=t^*(u(t,q),\pm q)$. 
The transformation can be found directly from $(S^*_{\a\b})$ in flat
coordinates by comparing the asymptotics of row sums in Theorem $5.2$ 
modulo $\h^{-2}$. 
Introduce $\dot{S}^*_{\a\b}(t^*,q)$ by
\[ S^*_{\a\b}=[\d_{\a\b}+\dot{S}^*_{\a\b}\h^{-1} +o(\h^{-1})] e^{t^*_{\b}/\h}
(e^0_{\b} e'_{\b})^{-1} (-1)^{\dim V} .\]
After some elementary computation we get
\[ t_{\b} = t_{\b}^* +\sum_{\a} e'_{\a} \dot{S}^*_{\a\b} (t^*,\pm q) 
(e'_{\b})^{-1} .\]
Notice that modulo $\h^{-2}$ the GW-potential $S^*_{\a\b}$ equals
\[ \d_{\a\b}(e^0_{\b}e'_{\b})^{-1} (-1)^{\dim V} + 
\h^{-1}  \sum_{d,n} q^d (\phi_{\a},t^*,...,t^*,\phi_{\b})/n! .\]
We obtain from this that in more invariant terms the change of variables is 
describes by the GW-invariants
\[ t_{\b}=\sum_{n,d} (\pm q)^d(e',t^*,...,t^*,\phi_{\b} e^0_{\b})^*/n! \]
where $e'=\sum_{\a} e'_{\a}\phi_{\a}$ is the equivariant Euler class of
$V$. In particular, the Jacobian of the change of variables is described
by the operator $e'\circ $ of quantum multiplication in 
$H^*_G(E^*)$:
\[ dt_{\b} = \sum_{\c} (e\circ )_{\c}^{\b} (t^*, \pm q) dt^*_{\c} .\]
We reiterate a question posed in \cite{Gi1}: how general is the 
nonlinear Serre duality relationship between GW-theory of dual supermanifolds
and bundle spaces?

\medskip

{\bf Toric supermanifolds.}
Let us assume now that $V$ is a direct sum of {\em positive} line bundles over
a toric symplectic manifold $X$ with equivariant Chern classes 
\[ v_j=\sum_{i=1} p_i l_{ij} -\l'_j ,\ j=1,...,l, \]
($p_i$ here are the same as in Theorem $4.2$)
and restrict ourselves to the study of
the fundamental solution $(S_{\a\b})$ for $\Pi E$ along $H^0\oplus H^2$
(with coordinates $t_0$ and $t=\log q$ respectively). As in Section $4$,
we have 
\[ \sum_{\a} S_{\a\b} =\lan J (q,\h) e^{(t_0+p\log q)/\h}, \phi_{\b} \ran \]
where $J$ is a formal vector $q$-series with coefficients which are rational
functions of $\l ,\l'$ and $\h$. Combining the nonlinear Serre duality  
with the mirror theorem for $V^*$ we conclude that the series $J$ is obtained
by a change of variables from its hypergeometric counterpart $I$. 
More precisely,
introduce the hypergeometric vector-function
\[ I = \sum_d q^d \Pi_{j=1}^l \frac{\Pi _{k=-\infty }^{L_j(d)} (v_j+k\h)}
{\Pi_{k=-\infty}^0 (v_j+k\h)} \Pi_{j=1}^m \frac{\Pi_{k=-\infty}^0 (w_j+k\h)}
{\Pi_{k=-\infty}^{D_j(d)} (w_j+k\h)} \]
whose terms differ from those in the hypergeometric series in Theorem $4.2$ 
(let us denote here that series by $I^*$)
by the factors $(-1)^l(v_1...v_l)^{-1}(v_1+L_1(d)\h)...(v_l+L_l(d)\h)$ and by
the signs $\pm q$. Taking into account Corollaries $4.3$ and $4.4$  
we arrive at the following mirror theorem for toric supermanifolds \cite{Gi3}.
\footnote{Our hypotheses here are somewhat more restrictive than in 
\cite{Gi3}, first --- because we assume that $V$ is {\em strictly} positive, 
and second --- because of the condition $\dim V >1$. In applications to 
hypersurfaces, say, in $\CC P^n$ the last condition is not constraining
because one can describe the same hypersurface as a codimension $2$ complete
intersection in $\CC P^{n+1}$. (We suggest the reader to consider the example
of Calabi-Yau $3$-folds given by two equations of degrees $5$ and $1$ in 
$\CC P^5$ in order to observe how the mirror transformation of \cite{COG} 
emerges from our formulas in the non-equivariant limit.)
We believe that the same trick can be applied to 
hypersurfaces in general toric varieties by extending GW-theory to 
K\"ahler orbifolds following M. Kontsevich's proposal.}  

\medskip

{\bf Corollary 5.3.} {\em Suppose that the toric supermanifold $\Pi E$ has 
non-negative $1$-st Chern class and that $\dim V>1$. Then the GW-potential
$J (q,\h) e^{(t_0+p\log q)/\h }$ is obtained from the hypergeometric vector
series $I (q,\h ) e^{(t_0+p\log q)/\h }$ by the division $I\mapsto I/\f (q)$
and by the change of variables 
\[ t_0 \mapsto t_0+\l g (q) +\l' g' (q) +f_0 , 
\ \log q_i \mapsto \log q_i + f_i (q) \] 
unambiguously determined by the asymptotics} 
\[ I=\f  [1+\h^{-1}(\l g +\l' g' +f_0+ f_1 p_1 +...+ f_r p_r) + 
o(\h^{-1}) ] .\]  

\medskip

{\em Proof.} The series $v_1...v_l Ie^{(p\log q)/\h}$ is actually obtained by 
differentiating $I^*e^{(p\log q)/\h}$, as in Corollary $4.4$,  
in the directions corresponding to the classes
\footnote{Notice that $v_j$ here correspond to $-v_j$ in Section $4$.} 
$-v_1,...,-v_l$ and changing $q$ to $\pm q$. According to Corollary $4.4$
the product $e'=v_1 \circ ...\circ v_l$ in the quantum cohomology algebra of
$E^*$ equals $v_1...v_l \f (\pm q)$. This identifies components of the ratio 
$\f ^{-1} I e^{(t_0+p\log q)/\h }$ with the GW-potentials
\[ \sum (\pm q)^d (e',\frac{e^{(t_0+p\log q)/\h}\phi_{\b} e^0_{\b}}{\h-c})^* \]
whose asymptotical terms of order $\h^{-1}$ determine the change
of variables prescribed by Theorem $5.2$.

\medskip

{\bf Elliptic GW-invariants.} Trying to extend GW-theory of convex 
sypermanifolds to higher genera in a way consistent with GW-theory for
corresponding complete intersections we encounter the following
difficulty. Let $f: (\S , \e)\to X$ be a stable map. Even if the bundle 
$V: E\to X$ is convex, the space $H^1(\S , f^*V)$ can be nontrivial depending
on the map. As a result, the spaces $H^0(\S , f^*V)$ do not form a vector 
bundle over $X_{g,n,d}$.  In any way, for a complete intersection $Y$ given
by a section of $V$ the virtual fundamental class $[Y_{g,n,d}]$ in $X_{g,n,d}$
does not have to be the cap-product of a cohomology class with the virtual
fundamental class $[X_{g,n,d}]$.

\medskip

{\em Counter-example.} Consider the moduli space $X_{1,1,0}=X\times \M_{1,1}$
of degree $0$ elliptic maps to an $m$-dimensional $X$. As we know, the 
virtual fundamental class is Poincar\'e-dual to 
$c_m(X)-\o \ c_{m-1}(X)$.
On the other hand the virtual fundamental class $[Y_{1,1,0}]$ for an 
(m-l)-dimensional submanifold
$i: Y \subset X$ is determined by the push-forwards $i_*(c_{m-l-1}(Y))$ and
$i_*(c_{m-l}(Y))$.
When $Y$ is given by a section of $V$, the push-forwards can be computed
in terms of $V$. Namely, the section identifies the normal bundle to $Y$ in 
$X$ with $i^*V$. Thus the (unstable) total Chern class of $Y$ equals
\[ \xsi^{m-l} + c_1(Y)\xsi^{m-l-1}+...+c_{m-l}(Y)=i^*
\frac{\xsi^m +c_1(X)\xsi^{m-1}+...+c_m(X)}{\xsi^l+c_1(V)\xsi^{l-1}+...+c_l(Y)}
 . \]    
The Chern classes of $Y$ we need are extracted from this ratio as
\[ c_{m-l}(Y)=i^*\Res_{\xsi=\infty} \frac{Chern (X)}{Chern (V)} 
\frac{d\xsi}{\xsi}\ \text{and}\ c_{m-l-1}(Y) =\Res_{\xsi=\infty }
\frac{Chern (X)}{ Chern (V)} \frac{d\xsi }{\xsi ^2} .\]
Since $\o ^2=0$, 
the virtual fundamental class $i_*(c_{m-l}(Y)-\o\ c_{m-l-1}(Y)$ of 
$[Y_{1,1,0}]$ is therefore computed as the Poincar\'e-dual to
\[  Euler (V) \Res_{\xsi =\infty} \frac{Chern (X)}{Chern (V)} 
\frac{d\xsi}{\xsi +\o} .\]
and has no reason to be a multiple of $c_m(X)-\o \ c_{m-1}(X)$.

\medskip

We believe that the virtual fundamental class $[Y_{g,n,d}]$ in
$X_{g,n,d}$ is the non-equivariant limit of a suitable equivariant homology 
class in $X_{g,n,d}$ and thus can be expressed via localization formulas in
terms of the fundamental class $[X^G_{g,n,d}]$ of the fixed point orbifold.
This would allow one to include complete intersections in the domain of 
applications of the results and methods of this paper. We are not ready
however to report upon any progress in this direction and hope to return to
this problem elsewhere.

\enddocument